\documentclass[english]{article}
\usepackage[T1]{fontenc}
\usepackage[latin9]{inputenc}
\usepackage[letterpaper]{geometry}
\usepackage{float}
\usepackage{amsmath}
\usepackage{graphicx}
\usepackage{color}
\usepackage{epstopdf}
\usepackage{adjustbox}
\usepackage{diagbox}
\usepackage{multirow}
\usepackage{array}
\usepackage{mathrsfs}
\usepackage{amsmath, amssymb, amsthm}
\usepackage{subfigure}
\usepackage{tikz}

\graphicspath{{./Figures/}}

\newtheorem{theorem}{Theorem}

\usepackage{amsfonts} 
\usepackage{color}
\usepackage{adjustbox}
\usepackage{diagbox}
\definecolor{black}{rgb}{0,0,0}

\definecolor{red}{rgb}{1,0,0}

\definecolor{blue}{rgb}{0,0,1}

\usepackage{multirow}

\makeatother

\usepackage{babel}
\usepackage{authblk}
\title{Constraint Energy Minimizing Generalized Multiscale Finite Element Method for high-contrast linear elasticity problem}

\author[1]{Shubin Fu \thanks{Corresponding Author}}
\author[1]{Eric T. Chung}
\affil[1]{Department of Mathematics, The Chinese University of Hong Kong, Hong Kong SAR}

\begin{document}
		\maketitle
\begin{abstract}
	In this paper, we consider the offline and online Constraint Energy Minimizing Generalized Multiscale Finite Element
	Method (CEM-GMsFEM) for high-contrast linear elasticity problem.
	Offline basis construction starts with
	an auxiliary multiscale space by solving local spectral problems. We select
	 eigenfunctions  that correspond to a few small eigenvalues to form the auxiliary space.  Using
	the auxiliary space, we solve a constraint energy minimization problem to construct offline multiscale spaces. The minimization problem is defined in
	the oversampling domain, which is larger than the target coarse block.
	To get a good approximation space, the oversampling domain should be 
	large enough. We also propose a relaxed minimization problem  to construct
	multiscale basis functions, which will yield more accurate and robust solution.  
	To take into account the influence of input parameters, such as source terms, we propose the construction of  online multiscale basis and an adaptive enrichment algorithm.
	We provide extensive numerical experiments on 2D and 3D models to show the performance of the proposed method. 
\end{abstract}	
\section{Introduction}
In many science and engineer problems, one encounters multiple scales and high
contrast. For example, wave propagation in fractured media, immiscible flow processes in poroelastic media and so on. Due to the advancement of
media characterization methods and geostatistical modeling techniques, 
the media can be  detailed at very fine scales, as a result, one needs to
solve huge dimensional algebraic systems. Therefore, model reduction methods 
are proposed by researchers to reduce the problem size and alleviate the computational  cost. Typical model reduction techniques include
 upscaling and multiscale methods. In upscaling methods \cite{wu2002analysis,durfolsky1991homo,gao2015numerical}, one typically upscales the media properties based on the homogenization theory so that the problem can be solved on a coarse grid.
 	In multiscale methods \cite{egw10,efendiev2009multiscale,chung2015mixed,chen2003mixed,jennylt03,Wheeler_mortar_MS_12, ArPeWY07,mortaroffline,hw97}, one still solves the problems on a coarse grid but with precomputed media dependent multiscale basis functions. 
 	
 	Among above mentioned multiscale methods, the multiscale finite element method (MsFEM) \cite{hw97,efendiev2009multiscale} is a classic multiscale method that has shown great success in various practical applications. However, the MsFEM assumes
that the media is scale separable. To overcome this assumption, the generalized 
 multiscale finite element method (GMsFEM)\cite{efendiev2013generalized} was proposed. The GMsFEM provide a   systematic way to construct multiple multiscale basis. In particular in GMsFEM, one first creates an appropriate snapshot space and then solve a carefully
designed local spectral problem in snapshot space. The basis space are filled with the dominant eigenvectors corresponding to small eigenvalues. 
 The GMsFEM's convergence depends on decay behavior of the eigenvalues
 of the local spectral problems \cite{efendiev2013generalized}.
In \cite{gmsfem_elasticity}, the authors applied the GMsFEM to solve the 
linear elasticity problem in high contrast problem, they consider both the 
continuous and discontinuous Galerkin method to couple the multiscale basis functions. In this paper, we will extend the recently proposed 
constraint energy minimizing GMsFEM (CEM-GMsFEM)\cite{chung2017constraint} for high contrast linear elasticity problem. 
The  CEM-GMsFEM consists of  two steps. One needs to first construct auxiliary basis functions by solving local spectral
problems. Then, for each auxiliary basis function, one can construct a multiscale basis via energy minimization problems on subdomains. We propose two versions . The first one is based on solving constraint energy minimization  problems and the second one is the relax version by solving
unconstrained energy minimization problems. The convergence of the CEM-GMsFEM not only depends on the eigenvalue but also depends on the coarse mesh size when the oversampling domain is carefully chosen.

To incorporate the influence of source and global media information, we also propose the construction of online multiscale basis.
The idea of online approach was first proposed in \cite{online_cg} and has been extended to various 
other cases (see \cite{online_dg,online_mixed,online_mortar}). The key idea
is using the residual information of the coarse-grid solution to construct multiscale basis.
These online multiscale basis functions can also be computed adaptively so that the error
can be decreased the most. The online basis of CEM-GMsFEM \cite{chung2018fast} will be computed in a 
oversampled domain, which is different from the original online approach \cite{online_cg}. 
We test our methods on 2D and 3D media with channels and inclusions.  By properly selecting the number
of basis functions and oversampling layers, we can observe that the multiscale solution can 
approximate the fine-scale solution accurately.

This paper will be organized as follows. In Section \ref{sec:pre}, we will  present some preliminaries.
In Section \ref{sec:offline}, the construction of offline multiscale basis functions of CEM-GMsFEM is discussed. In Section \ref{sec:online}, we present a online adaptive enrichment algorithm. In Section \ref{sec:convergence}, we provide some convergence results. In Section \ref{sec:numerical}, a few numerical results are
presented to demonstrate the performance of the method. Finally, some conclusions are given.
\section{Preliminaries}\label{sec:pre}
We consider isotropic linear elasticity problem in heterogeneous media as:
\begin{equation}
\label{ob_equ}
- \nabla \cdot  {\sigma(u)} =  {f}, \quad\text{ in } \; D
\end{equation}
where $D\subset\mathbb{R}^{d}(d=2,3)$ be a bounded domain representing the elastic body of interest, $u$ is the vector displacement field, $ {\sigma}( {u})$ is the stress tensor and it is related to the strain tensor $ {\epsilon}( {u})$ in the following way
\begin{equation*}
{\sigma(u)} = 2\mu  {\epsilon(u)} + \lambda \nabla\cdot  {u} \,  {\mathcal{I}},
\end{equation*}
where $\lambda>0$ and $\mu>0$ are the Lam\'e coefficients and can be highly
heterogeneous, $\mathcal{I}$ is  the identity tensor.
The strain tensor $ {\epsilon}( {u}) = (\epsilon_{ij}( {u}))_{1\leq i,j \leq d}$
is defined by
\begin{equation*}
{\epsilon}( {u}) = \frac{1}{2} ( \nabla  {u} + \nabla  {u}^T ),
\end{equation*}
where $\displaystyle \nabla  {u} = (\frac{\partial u_i}{\partial x_j})_{1\leq i,j \leq d}$.
In the component form, we have
\begin{equation*}
\epsilon_{ij}( {u}) = \frac{1}{2} \Big( \frac{\partial u_i}{\partial x_j} + \frac{\partial u_j}{\partial x_i} \Big), \quad 1\leq i,j \leq d. 
\end{equation*}
For simplicity, we will consider the homogeneous Dirichlet boundary condition $ {u} =  {0}$ on $\partial D$. 
Other types of boundary conditions can be taken care easily in the way used in classical approaches. 
 
\begin{figure}[ht]
	\centering
	\includegraphics[width=5in]{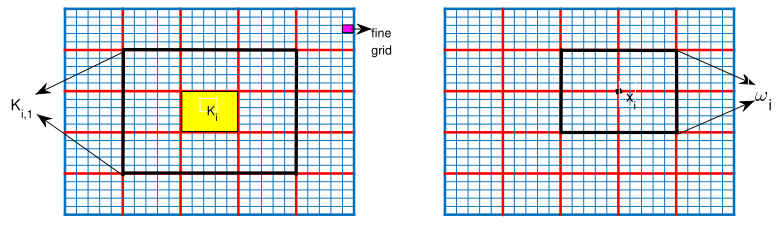}
	\caption{Illustration of the coarse grid, fine grid, oversampling domain $K_{i,1}$ and neighborhood.}
	\label{fig:grid}
\end{figure}

Let $\mathcal{T}^H$ be a conforming partition of the domain $D$. We call $\mathcal{T}^H$ the coarse grid and $H$ the coarse mesh size. 
Each element of $\mathcal{T}^H$ is called a coarse grid block. 
Denote $N_c$ be the total number of vertices of $\mathcal{T}^H$ and $N$ be the total
number of coarse blocks. Let $\{x_i\}_{i=1}^{Nc}$ be the set of vertices
in $\mathcal{T}^H$ and $\omega_i=\cup\{K_j\in\mathcal{T}^H|x_i\in\overline{K_j} \}$.
In addition, we let $\mathcal{T}^h$
be a conforming refinement of the triangulation $\mathcal{T}^H$.
We call $\mathcal{T}^h$ the fine grid and $h>0$ is the fine mesh size.
Figure \ref{fig:grid} shows an illustration of the fine scale grid, coarse scale grid, and oversampling domain.

Let $V=H_0^1(D)$, then the solution $u$ to (\ref{ob_equ}) satisfies
\begin{equation}
\label{cg_fine_sol}
a( {u},  {v}) = ( {f},  {v}), \quad \forall  {v}\in V,
\end{equation}
where
\begin{equation}
a( {u},  {v}) = \int_D \Big( 2\mu  {\epsilon}( {u}) :  {\epsilon}( {v}) 
+ \lambda \nabla\cdot  {u} \, \nabla\cdot  {v}  \Big) \; d {x},
\quad
( {f},  {v}) = \int_D  {f} \cdot  {v} \; d {x}
\end{equation}
and
\begin{equation}
{\epsilon}( {u}) :  {\epsilon}( {v}) = \sum_{i,j=1}^d \epsilon_{ij}( {u}) \epsilon_{ij}( {v}),
\quad
{f} \cdot  {v} = \sum_{i=1}^d f_i v_i.
\end{equation}
We will construct multiscale space $V_{ms}\subset V$ and obtain the 
solution in $V_{ms}$. That is, find $u_{ms}\in V_{ms}$ such that 
\begin{equation}\label{eq:ms}
a(u_{ms},v)=(f,v), \quad \forall  {v}\in V_{ms},
\end{equation}
To evaluate the accuracy of multiscale solution $u_{ms}$, we will compute
the solution of (\ref{cg_fine_sol}) on fine grid $\mathcal{T}_h$, denoted by
$u_h$ which is
fine enough to resolve all the heterogeneities of the exact $u$.
Let $V_h$ be the first-order Galerkin finite element basis space with respect to
the fine grid $\mathcal{T}_h$, and $q_1,\cdots,q_n$ be the basis set for 
$V_h$, then 
$u_h$ satisfies $A_hu_h=F_h$, where $A_h$ is a symmetric, positive definite matrix with
$A_{h,ij}=a(q_j,q_i)$, $F_h$ is a vector whose $i$-th component is $(f,q_i)$.
We will also use first order finite element on fine grid  $\mathcal{T}_h$ to compute the multiscale basis functions numerically. Then each multiscale 
basis function can be treated as a column vector $\Phi_i$, let $R=[\Phi_{1},\cdots,\Phi_{Nms}]$ be the 
matrix that stores all the multiscale basis functions (total number is $Nms$), then the  multiscale solution satisfies $u_{ms}=(R^TA_hR)^{-1} (R^TF_h)$, one can also project the coarse solution $u_{ms}$ into space $V_h$ by $u_{ms}^f=Ru_{ms}$.

\section{The construction of the CEM-GMsFEM basis functions.}\label{sec:offline}
This section is devoted to the construction of the multiscale basis functions.
There are two stages, the first stage is to construct the auxiliary multiscale basis function with the concept of generalized multiscale finite element method (GMsFEM). Then, we can construct the multiscale basis function by solving 
some energy minimizing problems in the oversampling domain.
\subsection{Auxiliary basis functions}
The auxiliary multiscale basis functions are constructed by solving a spectral problem in each coarse block $K_i$. More specifically, for each coarse block $K_i$, we let $V(K_i)$ be the restriction of $V$ on $K_i$, then we solve the spectral problem: find $(\lambda_j^i,\phi_j^i)\in \mathbb{R}\times V(K_i)$ 
\begin{equation*}
a_i(\phi_j^i,v)=\lambda_j^i s_i(\phi_j^i,v),\quad \forall v \in V(K_i)
\end{equation*}
where $a_i(u,v)=\int_{K_i}  2\mu  {\epsilon}( {u}) :  {\epsilon}( {v}) 
+ \lambda \nabla\cdot  {u} \, \nabla\cdot  {v}  $, $s_i(u,v)=\int_{K_i} \tilde{\kappa}uv$,
$\tilde\kappa = \sum_{i=1}^{N_c}(\lambda+2\mu) | \nabla \chi_i |^2$ and ${\chi_i}$ is a set of 
 partition of unity functions (see\cite{bm97})  on the coarse grid. 
 We choose  $J_i$ eigenfunctions  corresponding to first $J_i$ smallest eigenvalues  to form
 the local auxiliary space $V_{\text{aux}}(K_i)$, which is 
 \begin{equation*}
 V_{\text{aux}}(K_i) = \text{span} \{\phi_j^{i}|1\leq j\leq J_i \}.
 \end{equation*}
 We can normalize these eigenfunctions such that $s_i(\phi_j^i,\phi_j^i)=1$, we also denote $\Lambda$ as the 
 minimum first discarded eigenvalue. Then, the auxiliary space $V_{\text{aux}}$ is defined as the sum
 of all local auxiliary spaces $V_{\text{aux}}(K_i)$.
 We  define the notion of $\phi$-orthogonality in the space V. That is given a function 
 $\phi_j^i\in V_{\text{aux}}$, we say that a function $\psi\in V$ is $\phi_j^i$-orthogonal if 
 \begin{equation*}
s(\psi,\phi_j^i)=1,\quad s(\psi,\phi_{j^{\prime}}^{i^{\prime}})=0\quad \text{if }j^{\prime}\neq j\text{ or }
i^{\prime}\neq i.
 \end{equation*}
 where $s(u,v)=\sum_{i=1}^{N}s_i(u,v)$.
 We also define a projection operator $\pi$ from space $V$ to $V_{\text{aux}}$ by
 \begin{equation*}
 \pi(v)=\sum_{i=1}^{N}\sum_{j=1}^{J_i}s_i(v,\phi_j^i)\phi_j^i,\quad \forall v \in V,
 \end{equation*}

The kernel of the operator $\pi$ can be defined as 
 \begin{equation*}
 \tilde{V}=\{w\in V|\pi(w)=0\}.
 \end{equation*}
 \subsection{Offline multiscale basis functions }
 With the auxiliary space, we can introduce the construction of offline multiscale basis functions. For each coarse block $K_i$, we can extend
 this region by $m$ coarse grid layer and obtain an oversampled region 
 $K_{i,m}$ (see Figure \ref{fig:grid} for an example of $K_{i,1}$). Then for each auxiliary function $\phi_j^i\in V_{\text{aux}}$,
 the multiscale basis function $\psi_{j,ms}^i$ can be defined by 
 \begin{equation}\label{eq:uncon}
 \psi_{j,ms}^i=\text{argmin}\Big\{a(\psi,\psi)|\psi\in V_0(K_{i,m}), \psi\text{ is } \phi_{j}^i\text{-orthogonal} \Big\}.
 \end{equation}
 where $V_0(K_{i,m})=H_0^1(K_{i,m})$.
 By using Lagrange Multiplier, the problem (\ref{eq:uncon}) can be rewritten as the following
 problem: find $ \psi_{j,ms}^i\in V_0(K_{i,m}),$ $\lambda\in V_{\text{aux}}^i$ such that 
 \begin{equation}\label{eq:uncon1}
 \begin{split}
a(\psi_{j,ms}^i,p)+s(p,\lambda) &=0\quad \forall p\in V_0(K_{i,m}),\\
s(\psi_{j,ms}^i-\phi_j^i,q) &=0\quad\forall q\in V_{\text{aux}}^i(K_{i,m}),
 \end{split}
 \end{equation}
 where $V_{\text{aux}}^i(K_{i,m})$ is the union of all local auxiliary spaces for $K_j\subset K_{i,m}$.
 One can numerically solve above continuous problem with fine scale mesh.
 More specifically, denote $M_h$ be the matrix such that $M_{h,ij}=s(q_j,q_i)$,
 $A_h^i$and $M_h^i$ be the restriction of $A_h$ and $M_h$ on $K_{i,m}$ respectively.
 $P^i$ is the matrix that 
 includes all the discrete auxiliary basis in space $V_{\text{aux}}(K_{i,m})$.
 
 The matrix formulation of problem (\ref{eq:uncon1}) is 
\begin{equation}
\left[\begin{array}{cc}
A_h^i & M_h^{i}P^i \\
(M_h^{i}P^i)^T & 0\\
\end{array}\right]\left[\begin{array}{cc}\psi_h^i\\\lambda_{h}^i\end{array}\right] =\left[\begin{array}{cc}0\\I_i\end{array}\right]
\end{equation}
where $P^i_j$ is the $j$-th column of $P^i$, $\psi_{j,h}^i$ is discrete $\psi_{j,ms}^i$,
$I_i$ is a sparse matrix whose nonzero elements (all are 1) are in the diagonal of the matrix, the position of these nonzero elements depends on the index order of $K_i$ in $K_{i,m}$.

 Following \cite{chung2017constraint},  we can relax the $\phi$-orthogonality in (\ref{eq:uncon}) and get a relaxed version of the multiscale basis functions. More specifically, 
 we solve the following un-constrainted minimization problem:
 find $\psi_{j,ms}^i\in V_0(K_{i,m})$ such that
 \begin{equation}
 \psi_{j,ms}^i=\text{argmin}\Big\{a(\psi,\psi)+s(\pi\psi-\phi_j^i,\pi\psi-\phi_j^i)|\psi\in V_0(K_{i,m})\Big\}.
 \end{equation}
 which is equivalent to the following local problem
 \begin{equation}\label{eq:offline}
 a(\psi_{j,ms}^i,v)+s(\pi(\psi_{j,ms}^i),\pi(v))=s(\phi_j^i,\pi(v)),\quad 
 \forall v\in V_0(K_{i,m}).
 \end{equation}
Using above defined notation, 
then the matrix formulation of Equation (\ref{eq:offline}) is 
\begin{equation}
\left(A_h^i+M_h^{i}(P^iP^{i,T})M_h^{i,T}\right)\psi_{j,h}^i=P^i_jM_h^{i,T}
\end{equation}

\begin{figure}[H]
	\centering
	\subfigure[First mulitscale basis]{
		\includegraphics[width=2.5in]{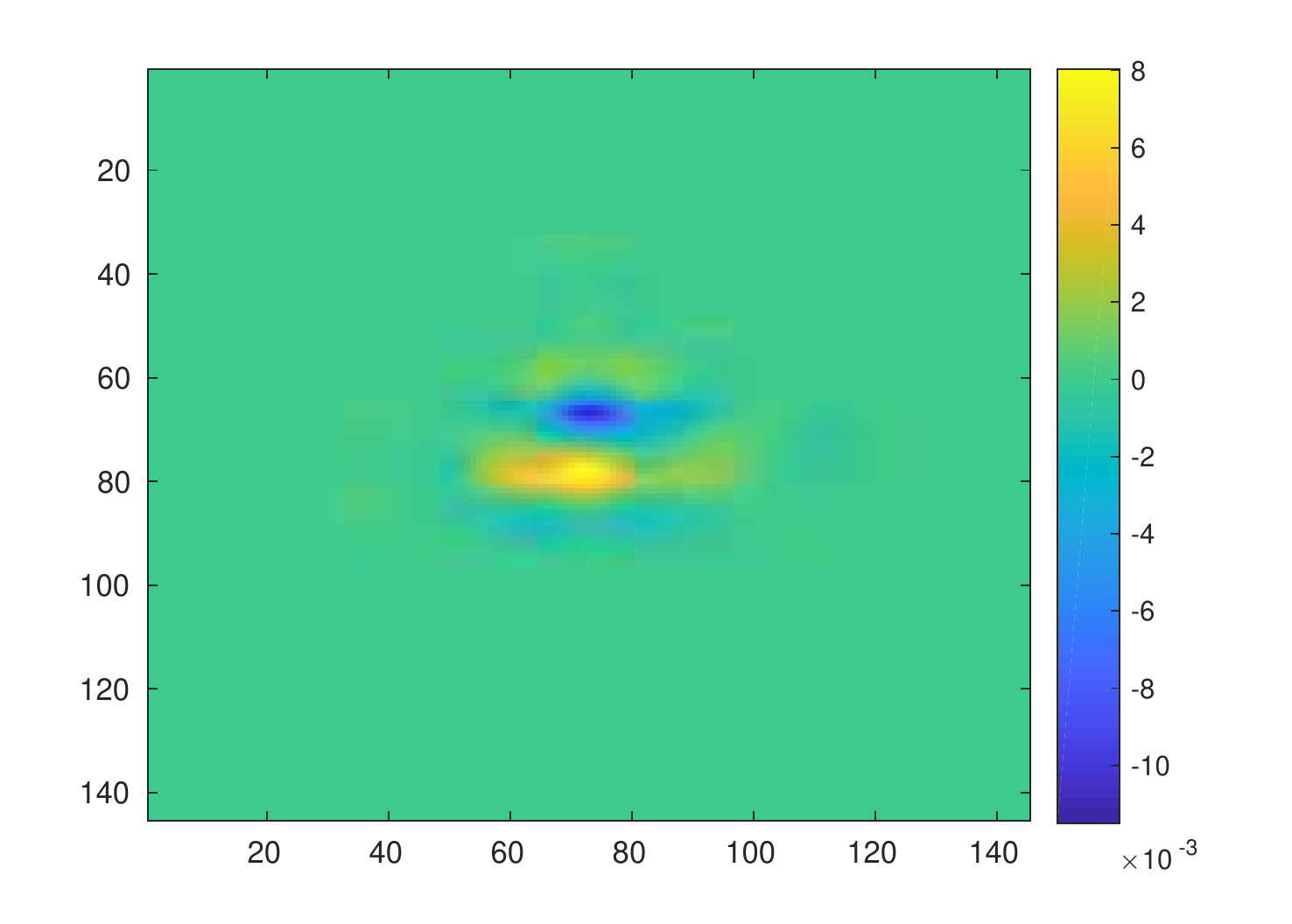}}
	\subfigure[Second mulitscale basis]{
		\includegraphics[width=2.4in]{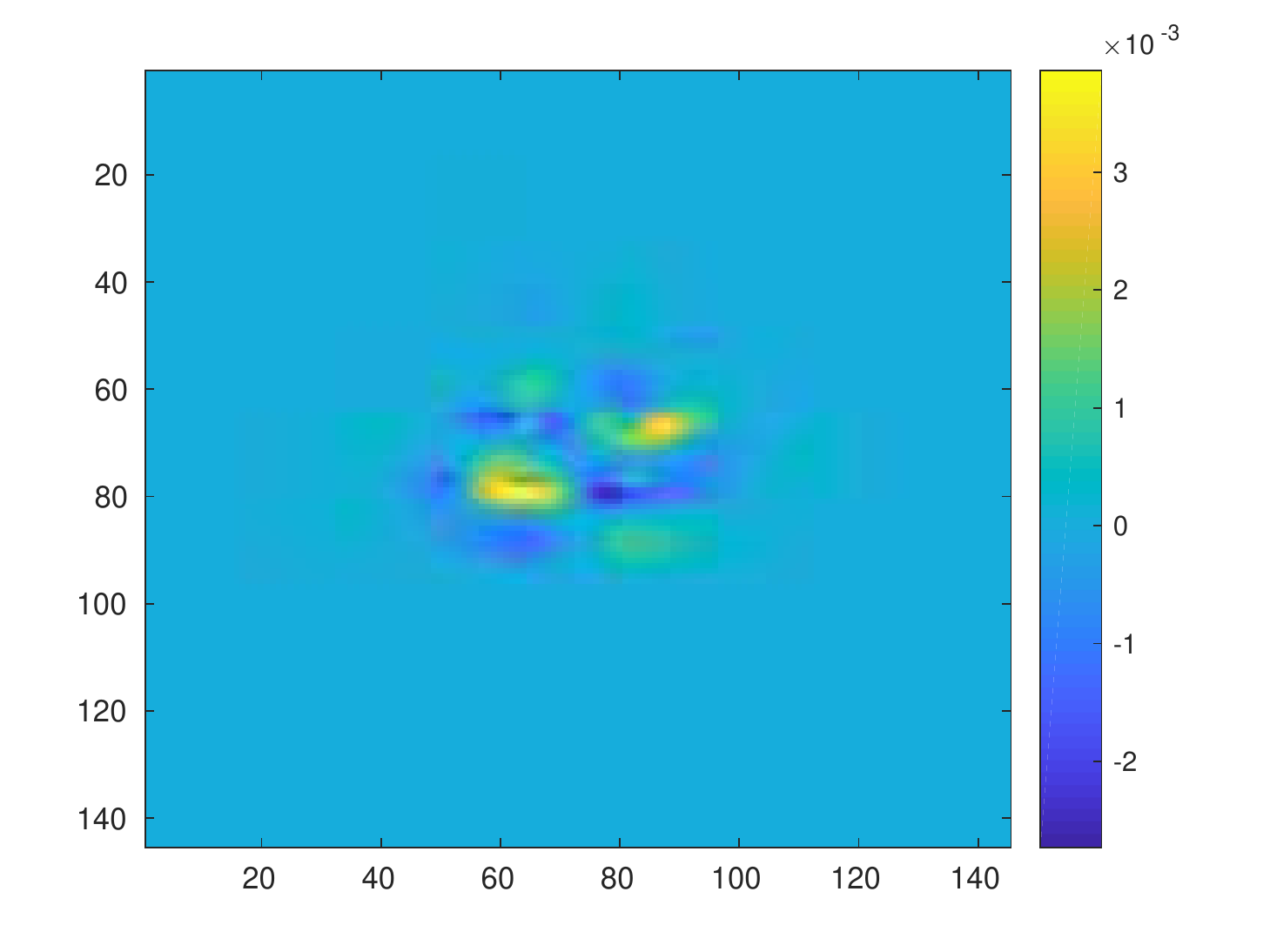}}
	\caption{Example of relaxed CEM-GMsFEM basis functions.}
	\label{fig:basis} 
\end{figure}

For each auxiliary basis $\phi_j^i$, one can get a multiscale basis $\psi^i_{j,ms}$. The final multiscale basis function space $V_{ms}$ is 
the span of all multiscale basis functions. Since the construction of the 
multiscale basis includes solving spectral problems and energy minimization
problems, therefore we call this method the CEM-GMsFEM. Figure \ref{fig:basis}
shows an example of relaxed CEM-GMsFEM basis functions, it can be observed that 
the multiscale basis concentrated on the support of auxiliary basis function
and decays outside the support. 
The multiscale basis functions can be treated as an approximation to global 
multiscale basis function $\psi_j^i\in V$ which is defined in a similar way, namely,
\begin{equation}
 \psi_j^i= \text{argmin}\Big\{a(\psi,\psi)|\psi\in V, \psi\text{ is } \phi_{j}^i\text{-orthogonal} \Big\}.
\end{equation} 
for the constraint case and
\begin{equation}
 \psi_j^i=\text{argmin}\Big\{a(\psi,\psi)+s(\pi\psi-\phi_j^i,\pi\psi-\phi_j^i)|\psi\in V\Big\}.
\end{equation}
for the relaxed case.
Then we can define the global space by $V_{glo}=\text{span}\{\psi_j^i\}$, this 
global basis functions decays exponential (see\cite{chung2017constraint}) and it is important to the convergence analysis.
\section{Online multiscale basis functions and adaptive enrichment}\label{sec:online}
In this subsection, we present the construction of online 
multiscale basis functions and an adaptive enrichment algorithm based on
an error estimate. Different with the offline basis, the online basis
functions are constructed iteratively using the residual of previous multiscale solution, therefore it contains the source information and global information of the media.

Let $u_{ms}\in V$ be the multiscale solution of (\ref{eq:ms}).
Then, we can define a residual functional $r:V\to\mathbb{R}$ by
\begin{equation}
r(v)=a(u_{ms},v)-\int_D fv,\quad \forall v\in V.
\end{equation}
The discrete residual in matrix form is $F_h-A_h\big(R((R^TA_hR)^{-1} (R^TF_h))\big)$.
For each neighborhood $\omega_i$ (see Figure \ref{fig:grid}), we can define the local residual functional $r_i:V\to \mathbb{R}$ by
\begin{equation}
r_i(v)=r(\chi_iv),\quad \forall v\in V.
\end{equation}.
The residual functional provides a way to measure the error $u-u_{ms}$ in $D$ and $\omega_i$.
Then, we can construction online basis function whose support is an oversampled region 
$\omega_i^+$ with the local residual $r_i$. More specifically, the online basis function
$\beta_{ms}^i$ satisfies following equation:
\begin{equation}\label{eq:online}
a(\beta_{ms}^i,v)+s(\pi(\beta_{ms}^i),\pi(v))=r_i(v),\quad \forall v\in H_0^1(\omega_i^+),
\end{equation}
Solving Equation (\ref{eq:online}) is similar with solving Equation (\ref{eq:offline}).
The online multiscale basis function is also localization results of corresponding global online basis
function $\beta_{glo}^i\in V$ defined by

\begin{equation}
a(\beta_{ms}^i,v)+s(\pi(\beta_{ms}^i),\pi(v))=r_i(v),\quad \forall v\in V.
\end{equation}
In practice, we can adaptively compute online basis for selected neighborhoods (with $i\in I$ for an
index set I).
After we construct the online basis functions, we can enrich the offline multiscale
basis space by adding the online basis, namely, $v_{ms}  =V_{ms}+\text{span}_{i\in I}\{\beta_{ms}^i\}$.
With the new multiscale basis function space, we can compute new multiscale solution and new 
basis space. These steps can be repeated until the residual norm is smaller than a given tolerance.  Before presenting the algorithm, we first define the $a$-norm $||\cdot||_a$ where $||u||_a^2=a( {u},  {u})$.
Next, we present the online adaptive enrichment algorithm.\\
\textbf{Online adaptive enrichment algorithm}

We first construct the offline basis functions space $V_{ms}^0$ introduced in Section \ref{sec:offline}.
We also choose a real parameter $\theta$ such that $0\leq\theta\leq 1$ to determine the
number of online basis functions added in each online iteration. Then for $m=0,1,2,\cdots$, we assume
that $V_{ms}^m$ is already obtained, then the updated multiscale basis functions space
$V_{ms}^{m+1}$.\\
Step 1: Find the multiscale in the current space $V_{ms}^m$. That is to find $u_{ms}^m\in V_{ms}^m$ such
that 
\begin{equation*}
a(u_{ms}^m,v)=(f,v),\quad \text{for all } v\in V_{ms}^m.
\end{equation*}\\
Step 2: For each neighborhood $\omega_i$, we compute the residual $z_i(v)$ by
\begin{equation*}
z_i(v)=a(u_{ms}^m,v)-(f,v),\quad \forall v\in V_0(\omega_i).
\end{equation*}
Denote $\delta_i=||z_i||_{a^*}=\text{sup}_{v\in V_0(\omega_i)}\frac{r(v)}{||v||_a}$. We rearrange the order 
of $\omega_i$ such that $\delta_1\geq\delta_2\cdots.$ Then we choose the first $k$ neighborhoods such that
\begin{equation*}
\sum_{i=k+1}^{N}\leq\theta_i^2\sum_{i=1}^{N}\delta_i^2
\end{equation*}
Step 3: Compute the local online basis functions in selected $k$ neighborhoods.
For each $1\leq i\leq k$ and neighborhood $\omega$, we find $\beta_{ms}^i\in V_0(\omega_i^+)$ satisfies
\begin{equation*}
a(\beta_{ms}^i,v)+s(\pi(\beta_{ms}^i),\pi(v))=r_i^m(v)\quad \forall v\in V_0(\omega_i^+),
\end{equation*}
where $r_i^m(v)=a(u_{ms}^m,\chi_iv)-\int_Df\chi_i v$\\
Step 4: Update the multiscale basis function space. That is form $V_{ms}^{m+1}$ by
\begin{equation*}
V_{ms}^{m+1}=V_{ms}^{m}+\text{span}_{1\leq i\leq k}{\beta_{ms}^i}.
\end{equation*}
\section{Convergence results}\label{sec:convergence}
In this section, we provide convergence results without giving the details of the analysis since it is quite similar with the techniques used in \cite{chung2018fast,chung2017constraint}. We define $s$-norm $||\cdot||_s$
by $ ||u||_s^2=\int_D \tilde{k}u^2$.
 We have following three theorems.
 \begin{theorem}
 	Let $u$ be the solution of equation (\ref{cg_fine_sol}) and $u_{ms}^{\text{off}}$ be the multiscale solution of (\ref{eq:ms}), the multiscale basis is the constraint case.   Then
 	we have
 	\begin{equation}
 	||u-u_{ms}^{\text{off}}||_a\leq C\Lambda^{-\frac{1}{2}}||\tilde{k}^{-\frac{1}{2}}f||_{L^2(D)}
 	+C(k+1)^{\frac{d}{2}}E^{\frac{1}{2}}||u_{glo}||_s
 	\end{equation}
 \end{theorem}
where $E$ is a constant that depends on $\Lambda$ and $m$, $u_{glo}$
 	 is the multiscale solution using corresponding global basis. 
\begin{theorem}
Let $u$ be the solution of equation (\ref{cg_fine_sol}) and $u_{ms}^{\text{off}}$ be the multiscale solution of (\ref{eq:ms}), the multiscale basis is the relaxed case. Then
we have
\begin{equation}
||u-u_{ms}^{\text{off}}||_a\leq C\Lambda^{-\frac{1}{2}}||\tilde{k}^{-\frac{1}{2}}f||_{L^2(D)}
+C(k+1)^{\frac{d}{2}}E^{\frac{1}{2}}(1+D)^{\frac{1}{2}}||u_{glo}||_s
\end{equation}
where $E$ and $D$ are constants that depend on $\Lambda$ and $m$, $u_{glo}$
is the multiscale solution using corresponding global basis. 
\end{theorem}
\begin{theorem}
Let $u$ be the solution of equation (\ref{cg_fine_sol}) and $u_{ms}^l$ be the sequence of multiscale solutions  generated by the 
online adaptive enrichment algorithm, the offline multiscale basis is the relaxed case Then
we have
\begin{equation*}
||u-u_{ms}^{l+1}||_a^2\leq 3(1+\Lambda^{-1})\big(C(m+1)^dE+2M^2\theta\big)
||u-u_{ms}^l||_a^2
\end{equation*}
where $E=3(1+\Lambda^{-1})(1+2(1+\Lambda^{\frac{1}{2}})^{-1})^{1-m}$,
$M$ is maximum number of overlapping subdomains and $C$ is a constant.
\end{theorem}
\section{Numerical results}\label{sec:numerical}
\begin{figure}[H]
	\centering
	\subfigure[Model 1]{
		\includegraphics[width=3in]{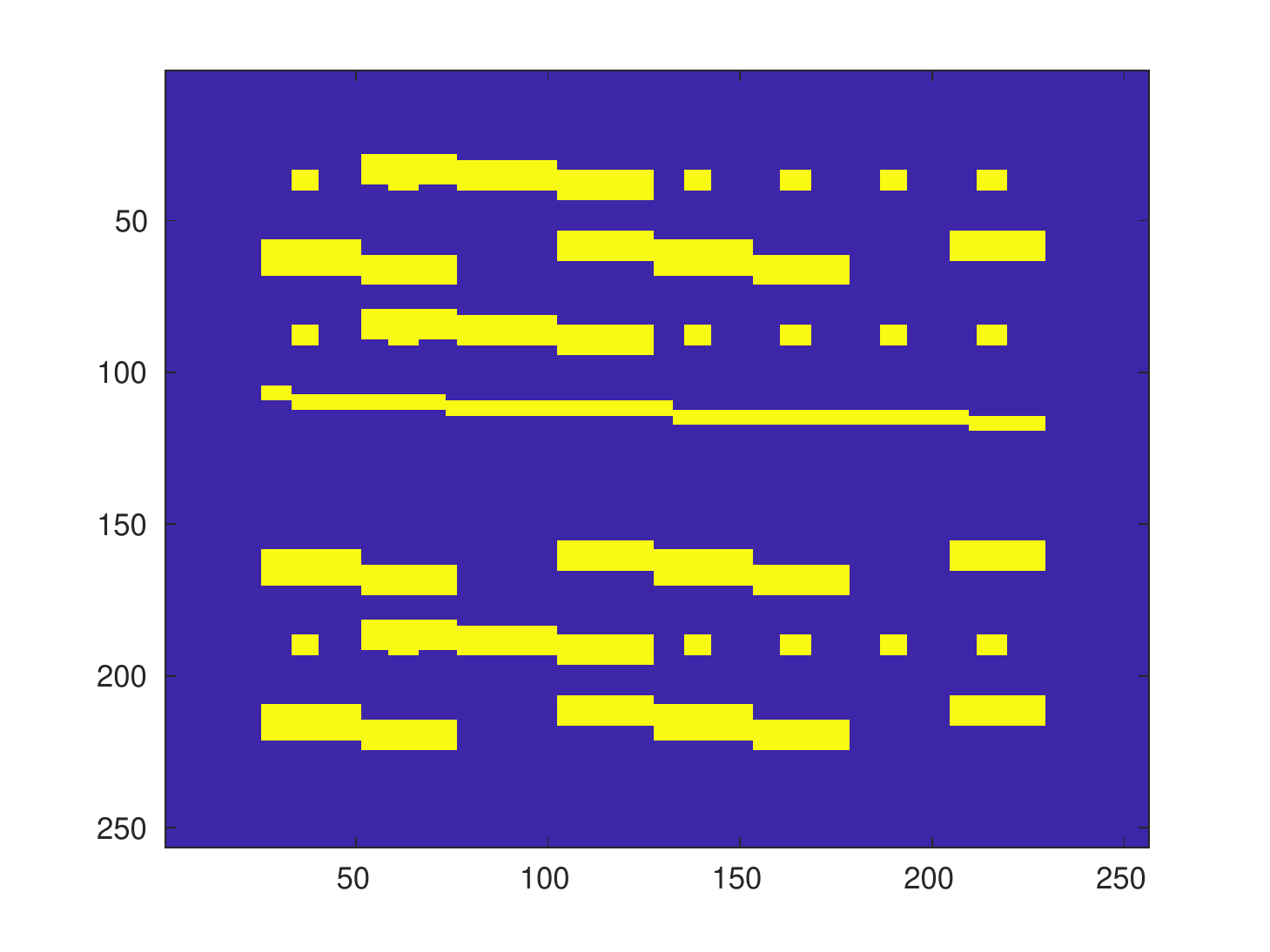}}
	\subfigure[Model 2]{
		\includegraphics[width=3in]{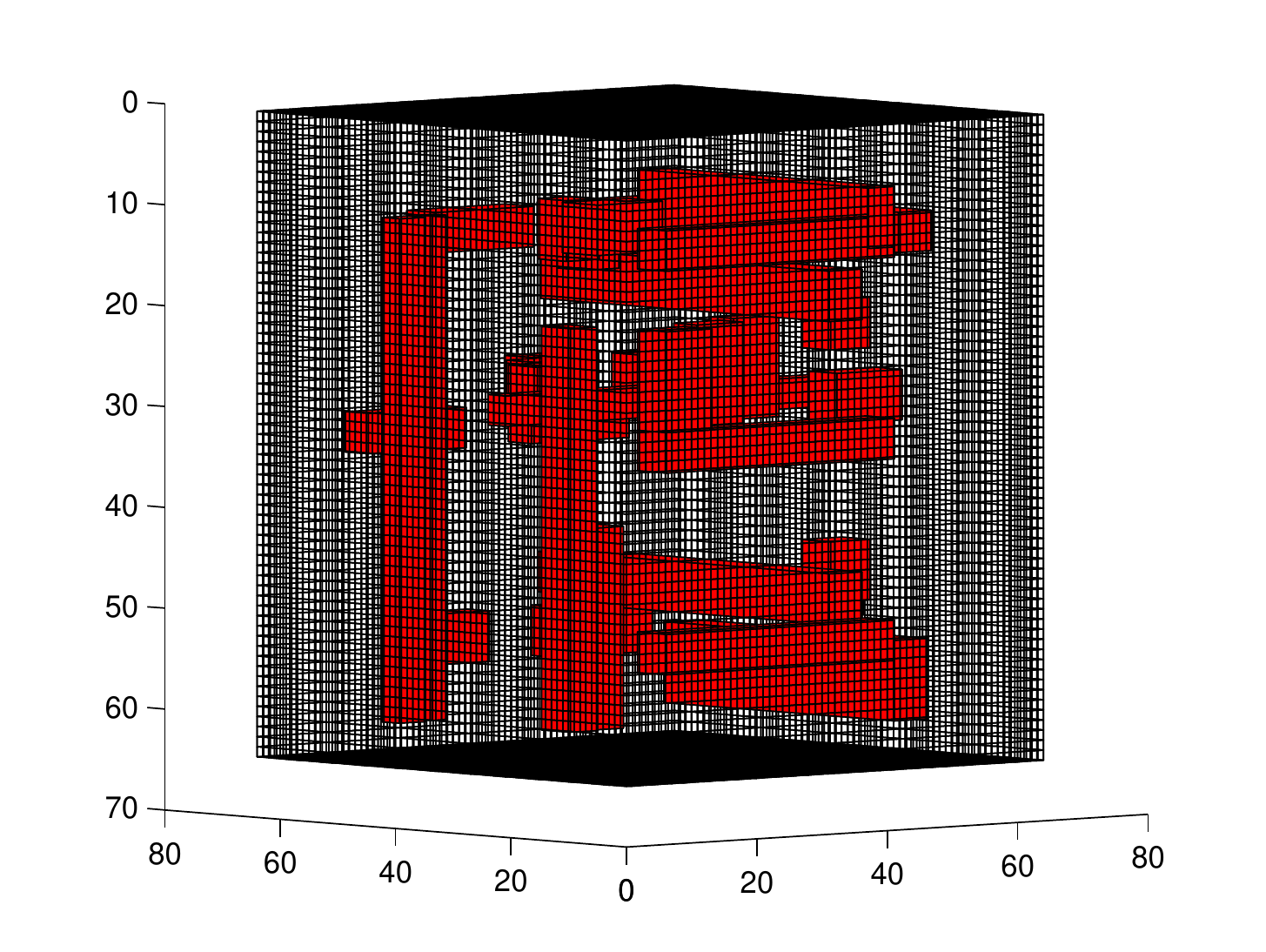}}
	\caption{Young's modulus.}
	\label{model} 
\end{figure}

In this section, we present several numerical experiments to show 
the performance of our method. The computational domain $D:= (0,1)^d$, we use constant force. 
We consider two high-contrast models whose 
Young's modulus $E(x)$ are depicted in Figure \ref{model}. As it is shown, both 
models contain  high conductivity channels and isolated inclusions.
We note that
for model 1, $E(x)=1$ in the blue region and $E(x)=E_1$ in the yellow region, while for model 2, $E(x)=1$ in the blank region and $E(x)=E_2$ in the red region. $\lambda(x)=\frac{\nu}{(1+2\nu)(1-\nu)}E(x), 
\mu(x)=\frac{1}{2(1+\nu)}E(x),$ the Poisson ration $\nu$ is $0.2$, both $E_1$ and $E_2$ equal $10^4$ unless
specifically illustrated.
The resolution of model 1 is $256\times256$, while for model 2 the resolution is $64\times64\times64.$
For all numerical results reported below, we 
use "$n_\text{ov}$" to represent the number of oversampling coarse layers used to compute the multiscale basis,
"$N_b$" is the number of basis used per coarse region, "Dof" means the degree of freedom of the resulting algebraic system, "$H$" is the coarse grid size. To quantify the accuracy of CEM-GMsFEM, we define 
 relative weighted $L^2$ norm error and weighted $H^1$ norm
 error as follows:
\begin{equation*}
e_{L^2}=\frac{||(\lambda+2\mu)(u_{ms}-u_h)||_{L^2(D)}}{||(\lambda+2\mu)u_h||_{L^2(D)}},\quad
e_{H^1}=\sqrt{\frac{a(u_{ms}-u_h,u_{ms}-u_h)}{a(u_h,u_h)}} 
\end{equation*}
where $u_h$ is the fine-grid first order FEM solution.
We first summarize our observations:
\begin{itemize}
	\item CEM-GMsFEM solution converges $H$ converges to 0 as $H$ converges to 0 for both relaxed and constraint case
	\item Relaxed CEM-GMsFEM is more accurate and robust than constraint CEM-GMsFEM under the same parameter setting
	\item Using more basis functions, adding mode oversampling coarse layers can improve the CEM-GMsFEM solution
	\item Online basis can accelerate the convergence of the CEM-GMsFEM solution.
\end{itemize}

\subsection{Constraint CEM-GMsFEM for model 1}
\begin{figure}[H]
	\centering
	\includegraphics[width=2.5in]{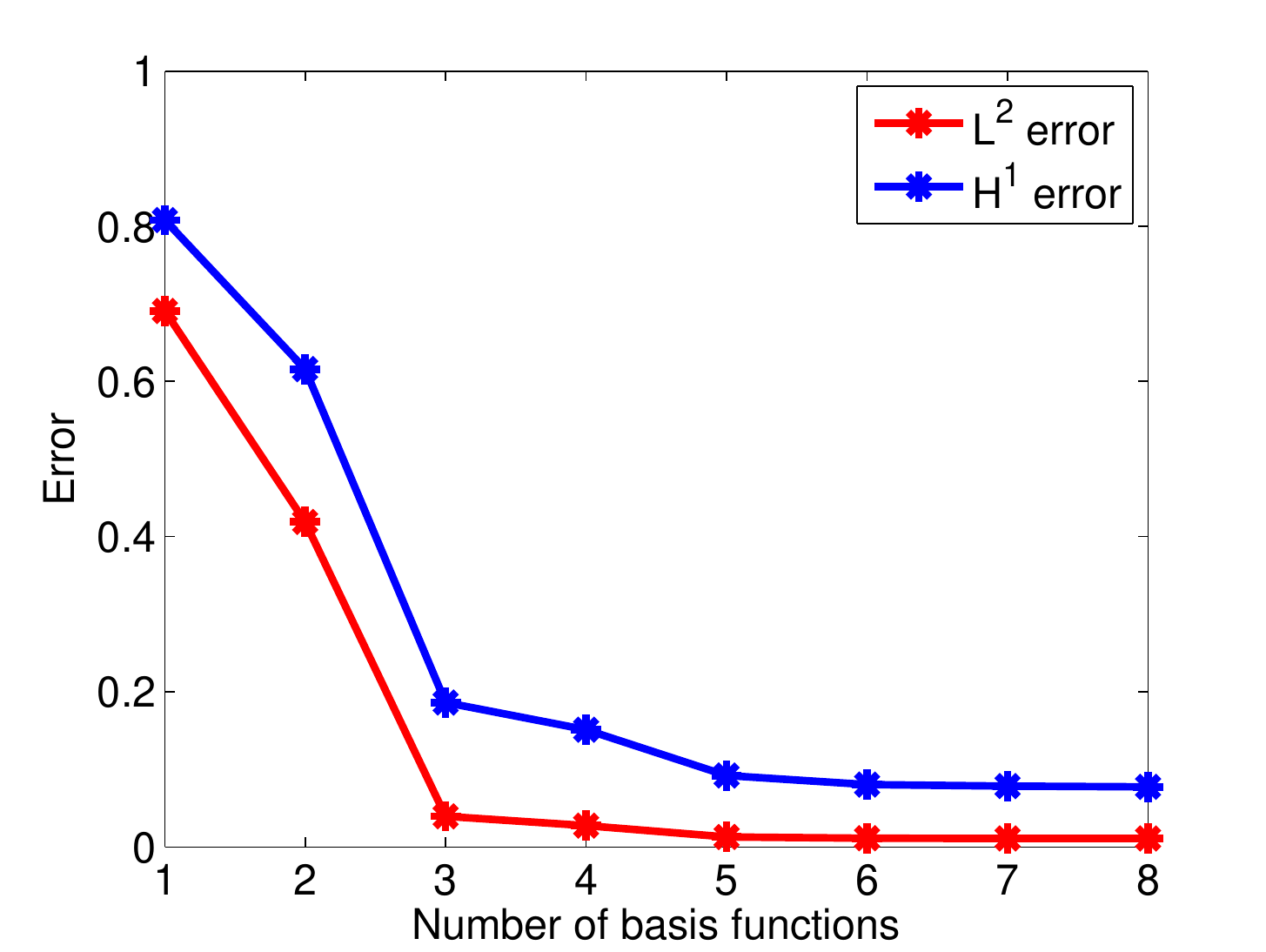}
	\caption{Numerical results  with different numbers of basis functions,  $H=1/16$, $n_\text{ov}=5$, constraint case.}
	\label{fig:convb}
\end{figure}

\begin{figure}[H]
	\centering
	\includegraphics[width=2.5in]{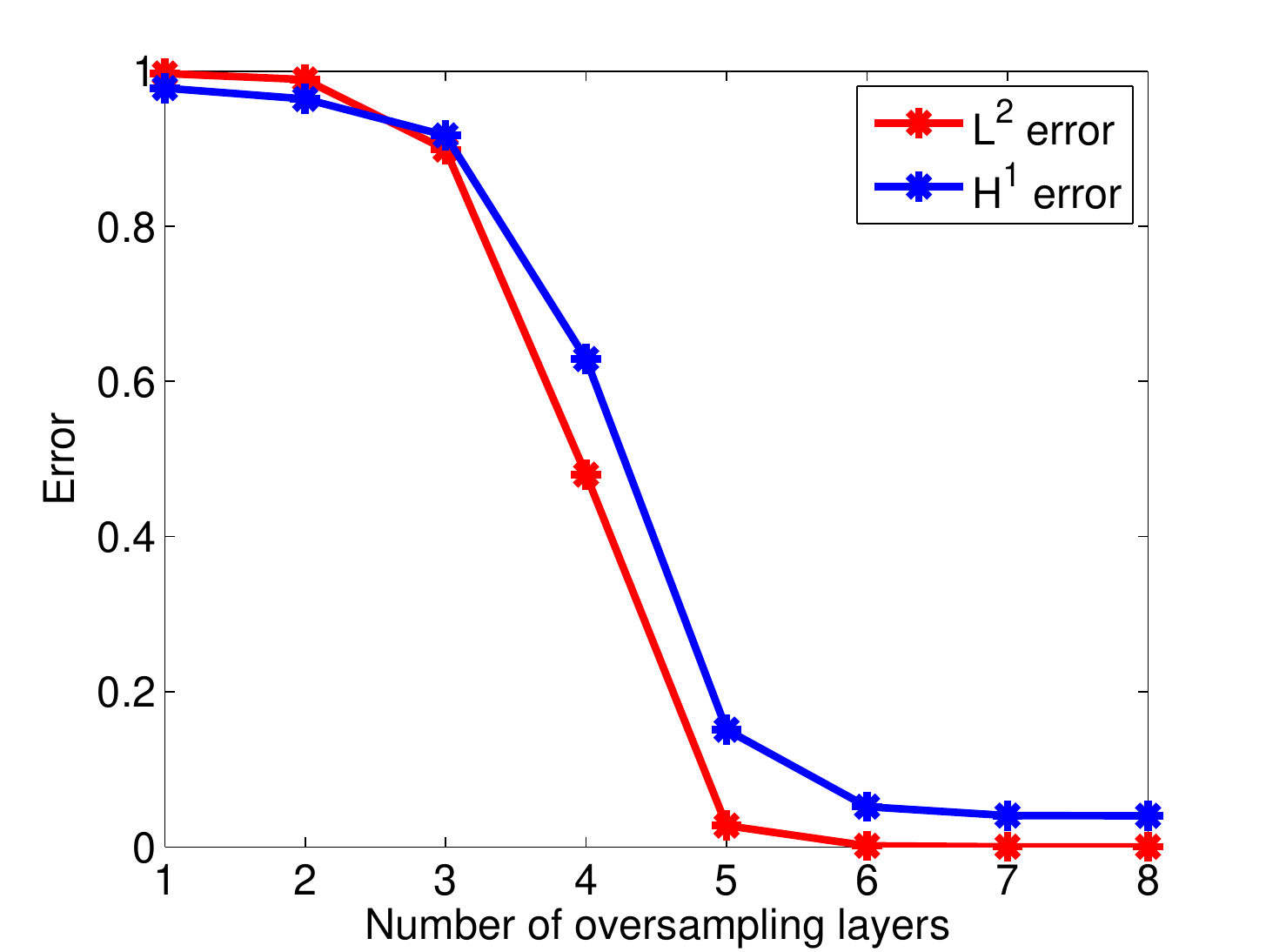}
	\caption{Numerical results with different numbers of oversampling layers, $H=1/16$, $N_b=4$, constraint case.}
	\label{fig:convov}
\end{figure}

We first present the test results of constraint CEM-GMsFEM on model 1. 
 The convergence history with various coarse mesh sizes $H$ are shown in Table \ref{ta:conH1} and Table \ref{ta:conH2}. For the simulation results reported in Table \ref{ta:conH1}, we take the number of oversampling layer to be approximately $3\text{log}(H)/\text{log}(1/8)$, as we can see, although the coarse solution converges as $H$ decreases, however the accuracy is not satisfiable, the $L^2$ error decrease from $68\%$ to only $26\%$. Therefore we increase the 
 number of oversampling layer to approximately $4\text{log}(H)/\text{log}(1/8)$, the corresponding results are reported in
 Table \ref{ta:conH2}. We find that both the $L^2$ and $H^1$ error improve a lot, and the errors decay also become faster.
  We emphasize that, in these simulations, we use 4 basis functions per coarse block since the
  eigenvalue problem on each coarse block has 4 small eigenvalues, and we need to include
  the first 4 eigenfunctions in the auxiliary space based on our theory.
  We also test the influence of number of basis, the results are presented in Figure \ref {fig:convb}.
 It can be seen clearly that increasing the number of basis will increase the accuracy of the CEM-GMsFEM solution.
 By varying the number of the oversampling coarse layers, we get results shown in Figure \ref{fig:convov}.
 We observe that the size of the subdomain to compute multiscale basis is quite important to the accuracy of CEM-GMsFEM,
 using more oversampling coarse layers will definite lead to more accurate coarse solution, 
 this agrees with the observations from Table \ref{ta:conH1} and \ref{ta:conH2}. However, after the 
 number of oversampling layers exceeds a certain number, the errors decay become slower.
  We also test different contrast case with fixed oversampling layer and number of basis functions, the relative $H^1$ error is shown
 in Table \ref{ta:concon}, we find that the performance of the scheme will deteriorate as the medium contrast increases, which is predicted by
 theoretical analysis. This motivates the propose of relaxed CEM-GMsFEM.
\begin{table}[H]
	\centering \begin{tabular}{|c|c|c|c|c|c|}\hline
		$N_b$ &$H$& $n_\text{ov}$  & $e_{L^2}$   & $e_{H^1}$   \tabularnewline\hline
		4&1/8	&3&6.77e-01&7.87e-01 \tabularnewline\hline
		4&1/16	&4&4.80e-01&6.28e-01 \tabularnewline\hline
		4&1/32	&5&3.70e-01&5.23e-01\tabularnewline\hline
		4&1/64	&6&2.62e-01&4.40e-01\tabularnewline\hline
	\end{tabular}
	\caption{Numerical results with varying coarse grid size $H$ for the test model 1, constraint case.}
	\label{ta:conH1}
\end{table}

\begin{table}[H]
	\centering \begin{tabular}{|c|c|c|c|c|c|}\hline
		$N_b$ &$H$&   $n_\text{ov}$  & $e_{L^2}$   & $e_{H^1}$   \tabularnewline\hline
		4&1/8	&4&5.65e-02&2.35e-01 \tabularnewline\hline
		4&1/16	&5&2.73e-02&1.50e-01 \tabularnewline\hline
		4&1/32	&7&4.47e-03&5.57e-02\tabularnewline\hline
		4&1/64	&8&2.67e-03&4.27e-02\tabularnewline\hline
	\end{tabular}
	\caption{Numerical results with varying coarse grid size $H$ for the test model 1, constraint case.}
	\label{ta:conH2}
	
\end{table}


\begin{table}[H] 
	\centering
	\begin{adjustbox}{max width=\textwidth}	
		\begin{tabular}{|c|c|c|c|c|c|c|c}\hline
			\diagbox{$n_\text{ov}$}{$E_1$} &$10^2$ &$10^4$ &$10^6$&$10^8$&$10^{10}$  \tabularnewline\hline
			5&1.90e-02&8.04e-02&4.04e-01&6.94e-01&7.18e-01\tabularnewline\hline
			6&1.47e-02&2.38e-02&1.33e-01&5.19e-01&7.10e-01\tabularnewline\hline
			7& 1.45e-02&1.93e-02&3.43e-02&2.31e-01&6.07e-01\tabularnewline\hline
		\end{tabular}
	\end{adjustbox}
	\caption{Comparison ($e_{H^1}$) of various number of oversampling layers and different contrast values for test model 1,  constraint case, $H=1/16, N_b=6$.}
\label{ta:concon}
\end{table}

\subsection{Relaxed CEM-GMsFEM for model 1}
\begin{figure}[H]
	\centering
	\subfigure[First component]{
		\includegraphics[width=2.5in]{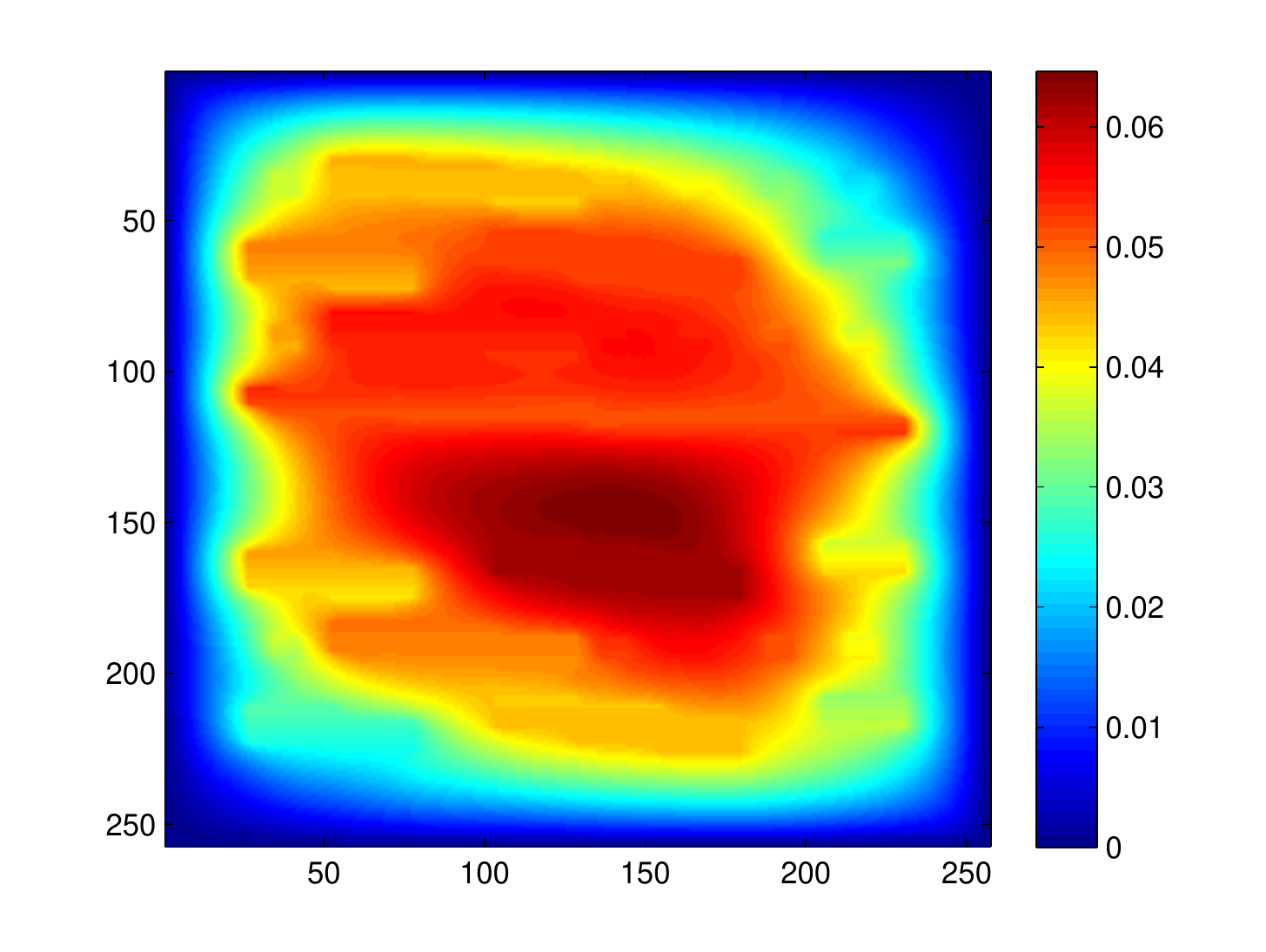}}
	\subfigure[Second component]{
		\includegraphics[width=2.5in]{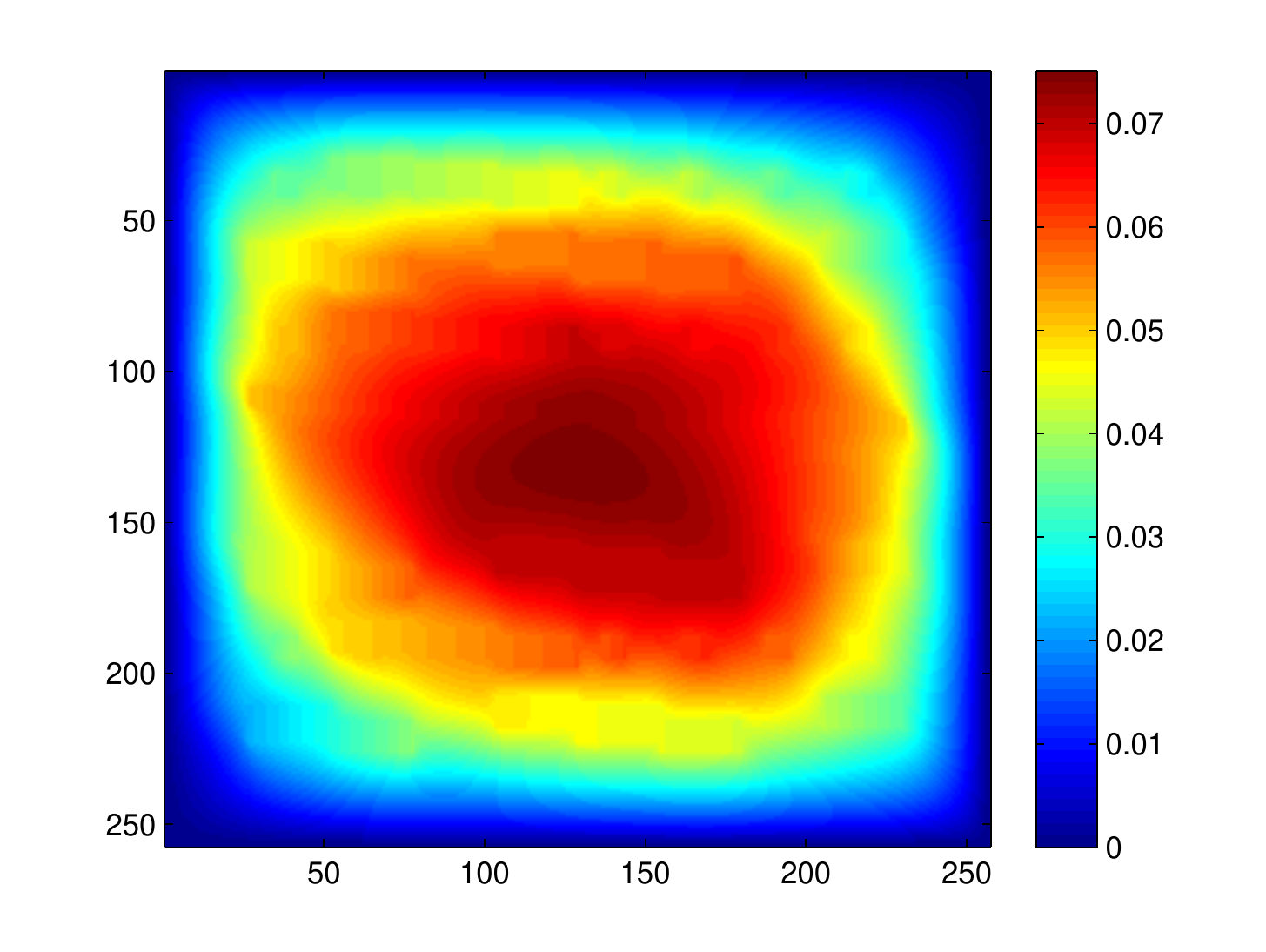}}
	\caption{Reference solution}
	\label{fig:finedis} 
\end{figure}

\begin{figure}[H]
	\centering
	\subfigure[First component]{
		\includegraphics[width=2.5in]{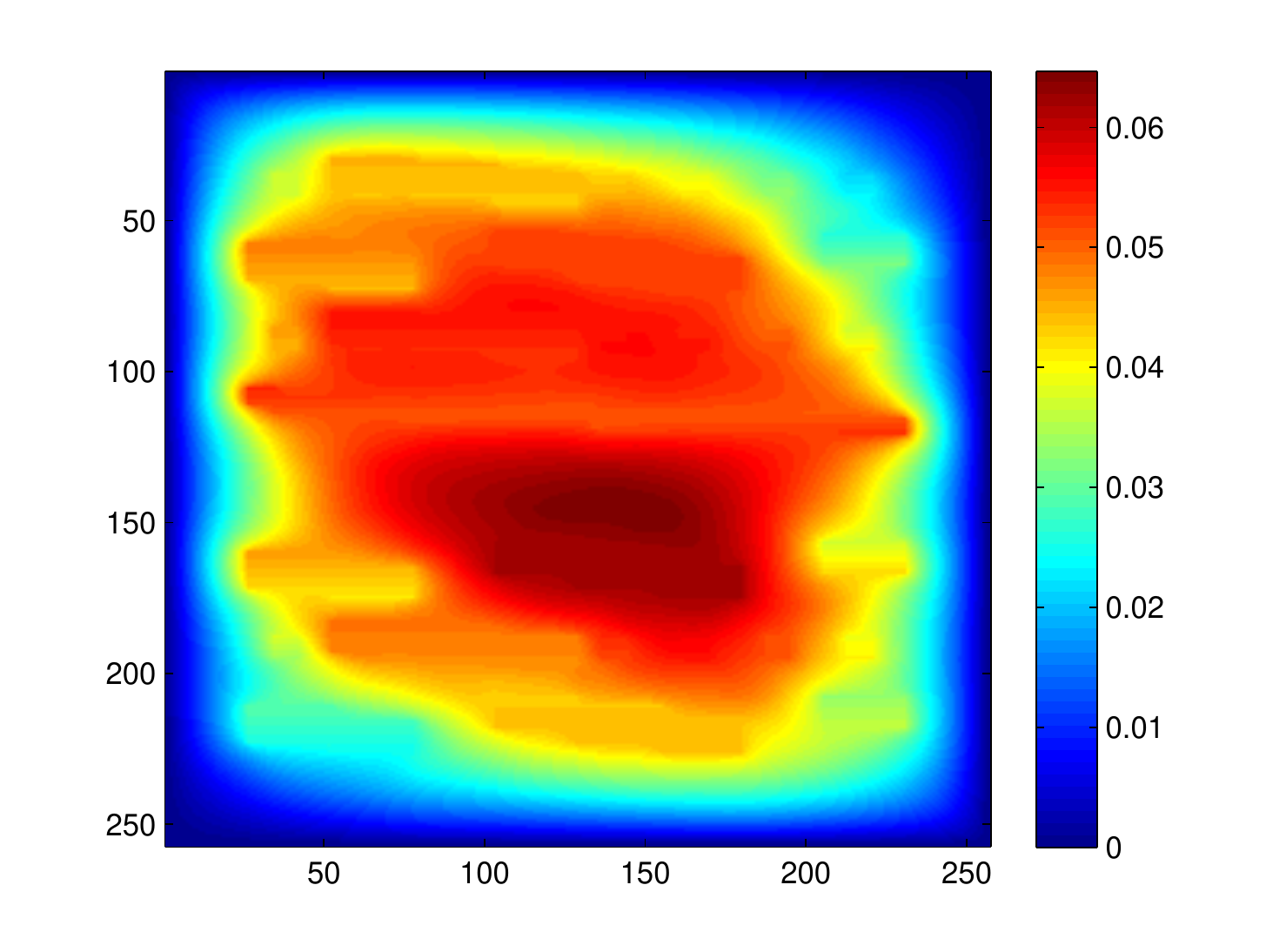}}
	\subfigure[Second component]{
		\includegraphics[width=2.5in]{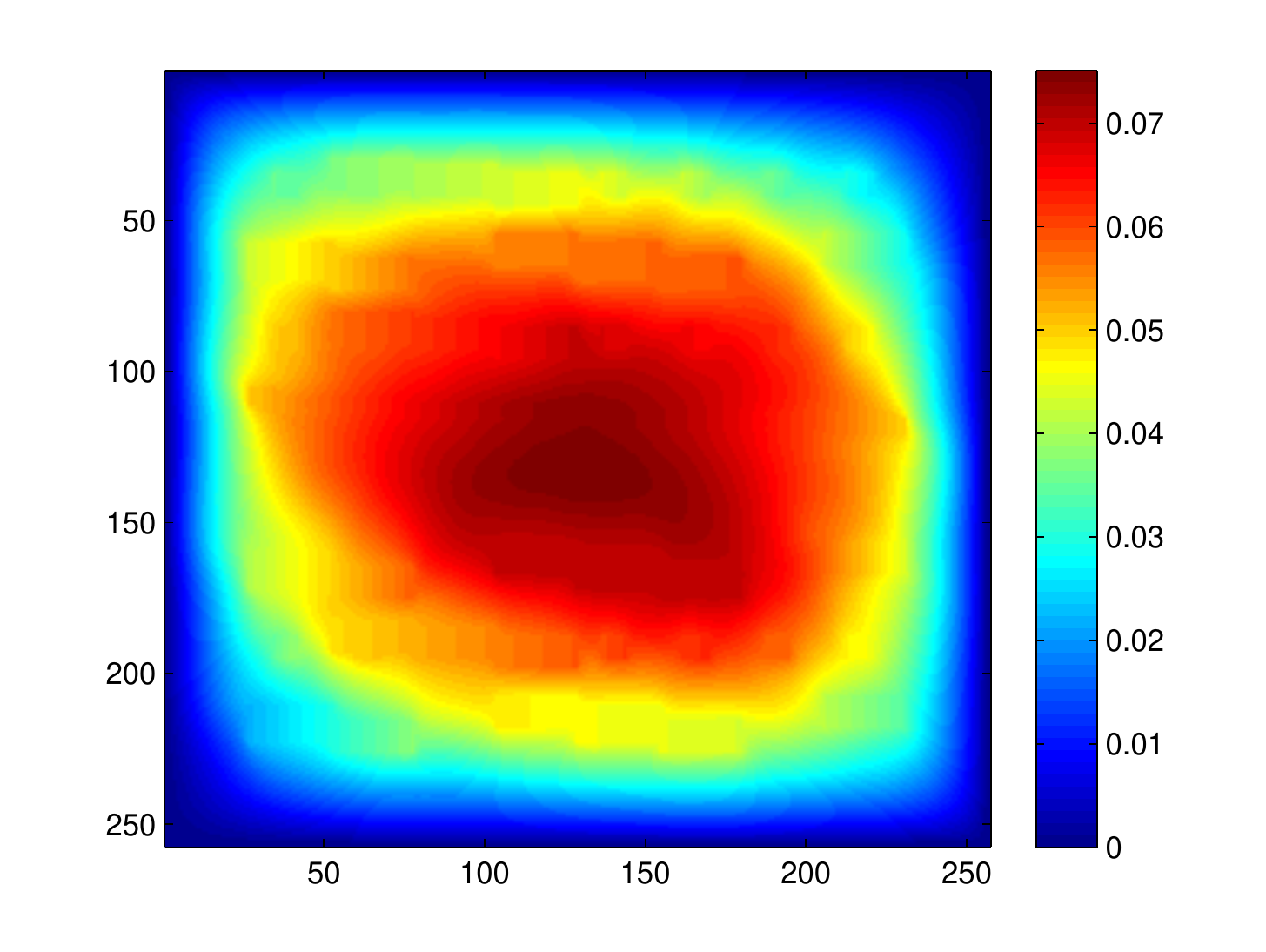}}
	\caption{Relaxed CEM-GMsFEM solution, $N_b$=4, H$=1/16$, $n_\text{ov}=4$.}
	\label{fig:msdis} 
\end{figure}

\begin{figure}[H]
	\centering
	\includegraphics[width=2.5in]{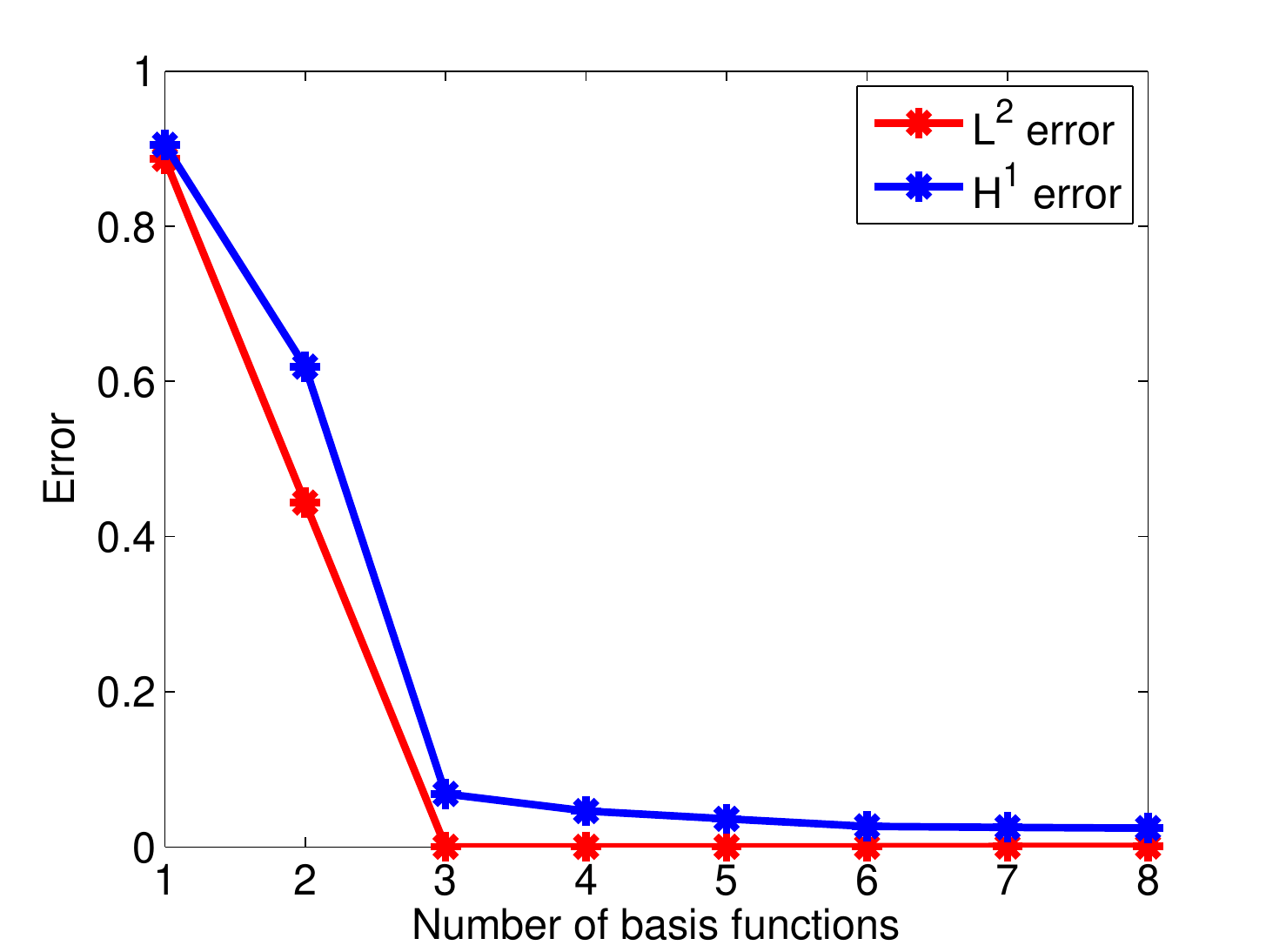}
	\caption{Numerical results  with different numbers of basis functions,  $H=1/16$, $n_\text{ov}=4$, relaxed case.}
	\label{fig:relaxvb}
\end{figure}

\begin{figure}[H]
	\centering
	\includegraphics[width=2.5in]{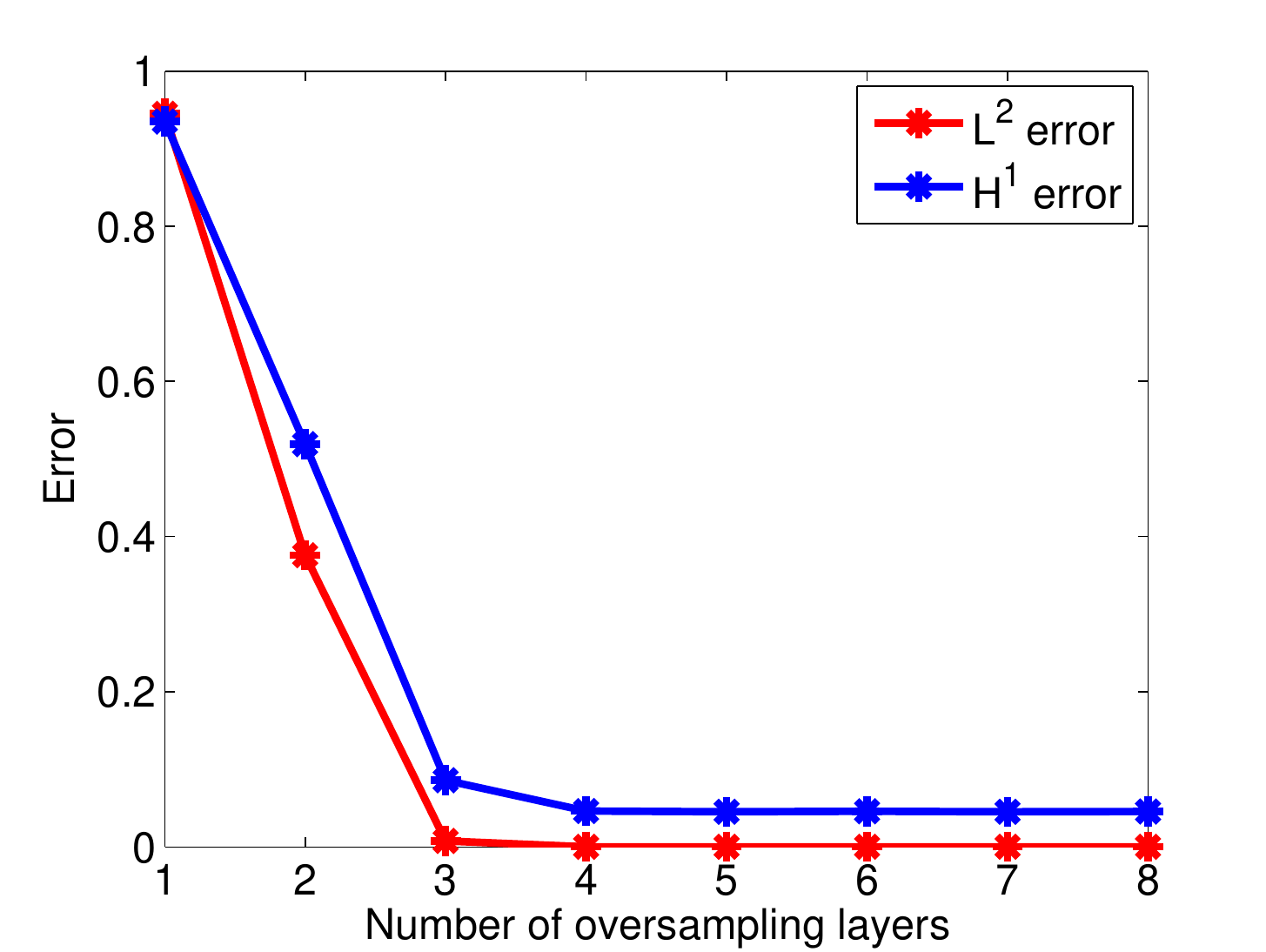}
	\caption{Numerical results with different numbers of oversampling layers, $H=1/16$, $N_b=4$, relaxed case.}
	\label{fig:relaxvov}
\end{figure}

In this subsection, we present the performance of Relaxed CEM-GMsFEM for 
model 1. We first linearly decrease the coarse-grid size and the results 
are shown in Table \ref{ta:relaxH1}. We observe that the coarse solution
converge fast as $H$ decreases, for example , the $L^2$ error decays from
$9\%$  to $0.0006\%$. The $L^2$ error convergence faster (close to 
second order) than the energy error (close to first order).
By comparing to the
similar test case in Table \ref{ta:conH1} and \ref{ta:conH2}, we see that the relaxed version needs fewer oversampling layers and obtains much better results.
Figure \ref{fig:finedis} shows the displacement fields of the reference solution, we can see complicated multiscale behavior of the solution. Figure \ref{fig:msdis} is the CEM-GMsFEM solution, we see the coarse solution can
capture almost all the details of the reference solution and there is almost 
no difference with the reference solution. 
We also investigate the performance with different number of eigenfunctions in the auxiliary space and number of oversampling layers, the results are reported 
in Figure \ref{fig:relaxvb} and Figure \ref{fig:relaxvov}. Again, as predicted by the theory, using more basis and larger subregion size will improve the accuracy of the CEM-GMsFEM solution.
The results of robustness test are shown in Table \ref{ta:relaxcon}, we can see that the relaxed CEM-GMsFEM is more robust with respect to the contrast.

\begin{table}[H]
	\centering \begin{tabular}{|c|c|c|c|c|c|}\hline
		$N_b$ &$H$&   $n_\text{ov}$  & $e_{L^2}$   & $e_{H^1}$   \tabularnewline\hline
		4&1/8	&3&9.04e-02& 2.72e-01 \tabularnewline\hline
		4&1/16	&4&1.35e-04& 4.63e-02 \tabularnewline\hline
		4&1/32	&5&2.20e-05& 1.47e-02\tabularnewline\hline
		4&1/64	&6&5.92e-06& 4.01e-03\tabularnewline\hline
	\end{tabular}
	\caption{Numerical results with varying coarse grid size $H$ for the test model 1, relaxed case.}
	\label{ta:relaxH1}
\end{table}


\begin{table}[H] 
\centering
\begin{adjustbox}{max width=\textwidth}	
\begin{tabular}{|c|c|c|c|c|c|c|c}\hline
\diagbox{$n_\text{ov}$}{$E_1$} &$10^2$ &$10^4$ &$10^6$&$10^8$&$10^{10}$  \tabularnewline\hline
5&1.98e-02&2.25e-02&2.36e-02&4.59e-02&3.00e-01\tabularnewline\hline
6&2.00e-02&2.26e-02&2.35e-02&2.33e-02&4.32e-02\tabularnewline\hline
7&1.97e-02&2.25e-02&2.26e-02&2.34e-02&2.35e-02\tabularnewline\hline	
\end{tabular}
\end{adjustbox}
\caption{Comparison ($e_{H^1}$) of various number of oversampling layers and different contrast values for test model 1 , constraint case, relaxed case, $H=1/16, N_b=6$.}
\label{ta:relaxcon} 
\end{table}

 We also test the online iterative algorithm, the results are reported in 
 Table \ref{ta:online1} and Table \ref{ta:online2}. We can see that with online basis
 functions, the convergence is very fast.
 Table \ref{ta:online1}
 shows the results of uniform enrichment, by comparing it with Table \ref{ta:online2}, we conclude that adaptive enrichment is better 
 especially in reducing the $L^2$ error.

\begin{table}[H]
	\centering \begin{tabular}{|c|c|c|c|c|c|}\hline
		Dof&  $n_\text{ov}$  & $e_{L^2}$   & $e_{H^1}$   \tabularnewline\hline
		16384	&4&1.55e-02&1.02e-01 \tabularnewline\hline
		20353	&4&2.91e-05&5.10e-03   \tabularnewline\hline
		24322	&4&3.54e-07&3.88e-04   \tabularnewline\hline
	\end{tabular}
	\caption{Uniform enrichment error decay history for the test model 1, $H=1/64$, 4 offline basis used.}
	\label{ta:online1} 
\end{table}

\begin{table}[H]
	\centering \begin{tabular}{|c|c|c|c|c|c|}\hline
		Dof&  $n_\text{ov}$  & $e_{L^2}$   & $e_{H^1}$   \tabularnewline\hline
		16384	&4&1.55e-02&1.02e-01 \tabularnewline\hline
		17112&4&9.25e-05&9.11e-03  \tabularnewline\hline
		17837&4&1.29e-05&3.91e-03    \tabularnewline\hline
		18574&4&1.02e-05&3.45e-03   \tabularnewline\hline
	\end{tabular}
	\caption{Adaptive enrichment with $\theta=0.1$ error decay history for the test model 1, $H=1/64$, 4 offline basis used.}
	\label{ta:online2} 
\end{table}

\subsection{Relaxed CEM-GMsFEM for model 2}
\begin{figure}[H]
	\centering
	\subfigure[First component]{
		\includegraphics[width=2in]{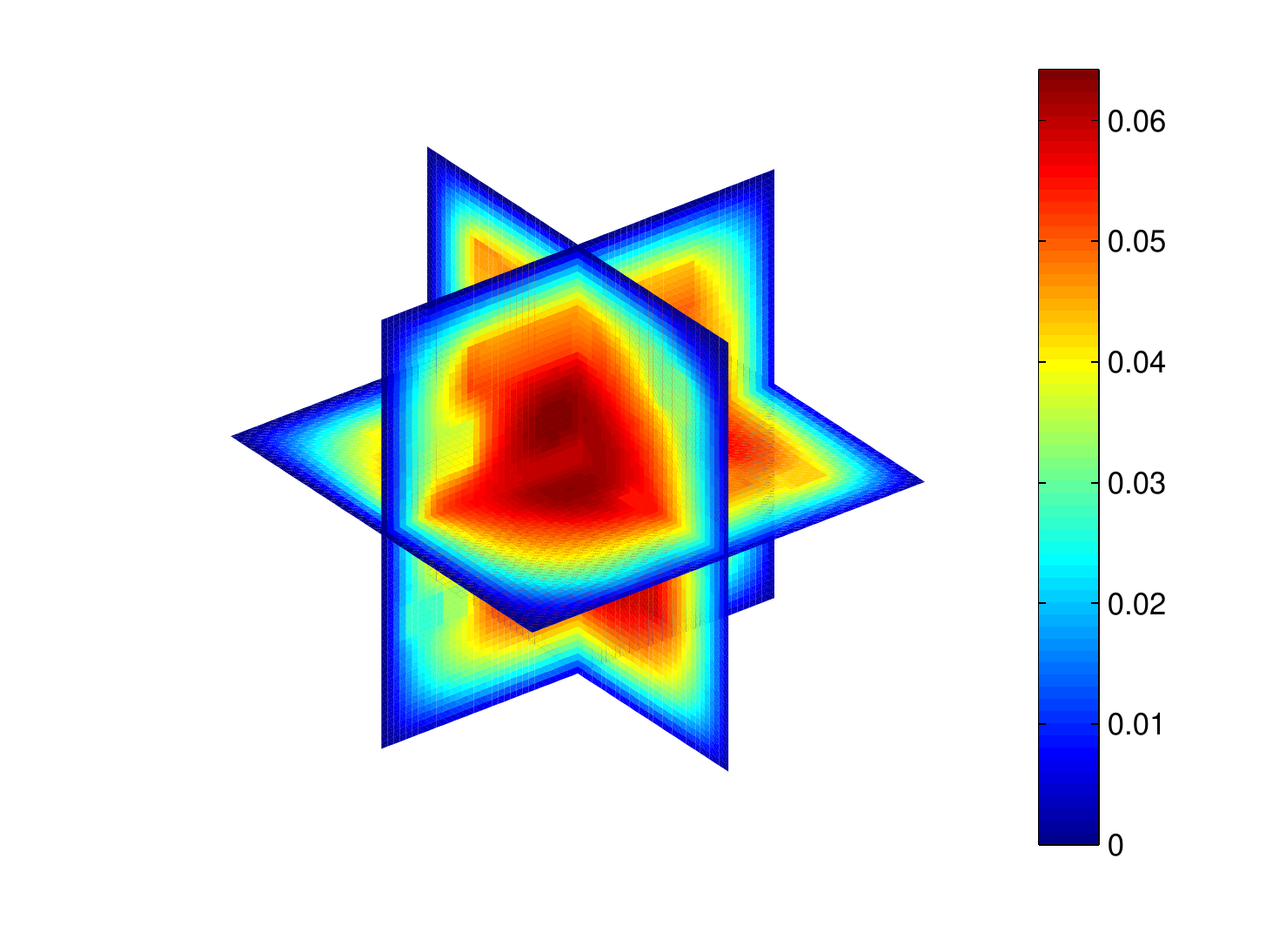}}
	\subfigure[First component]{
		\includegraphics[width=2in]{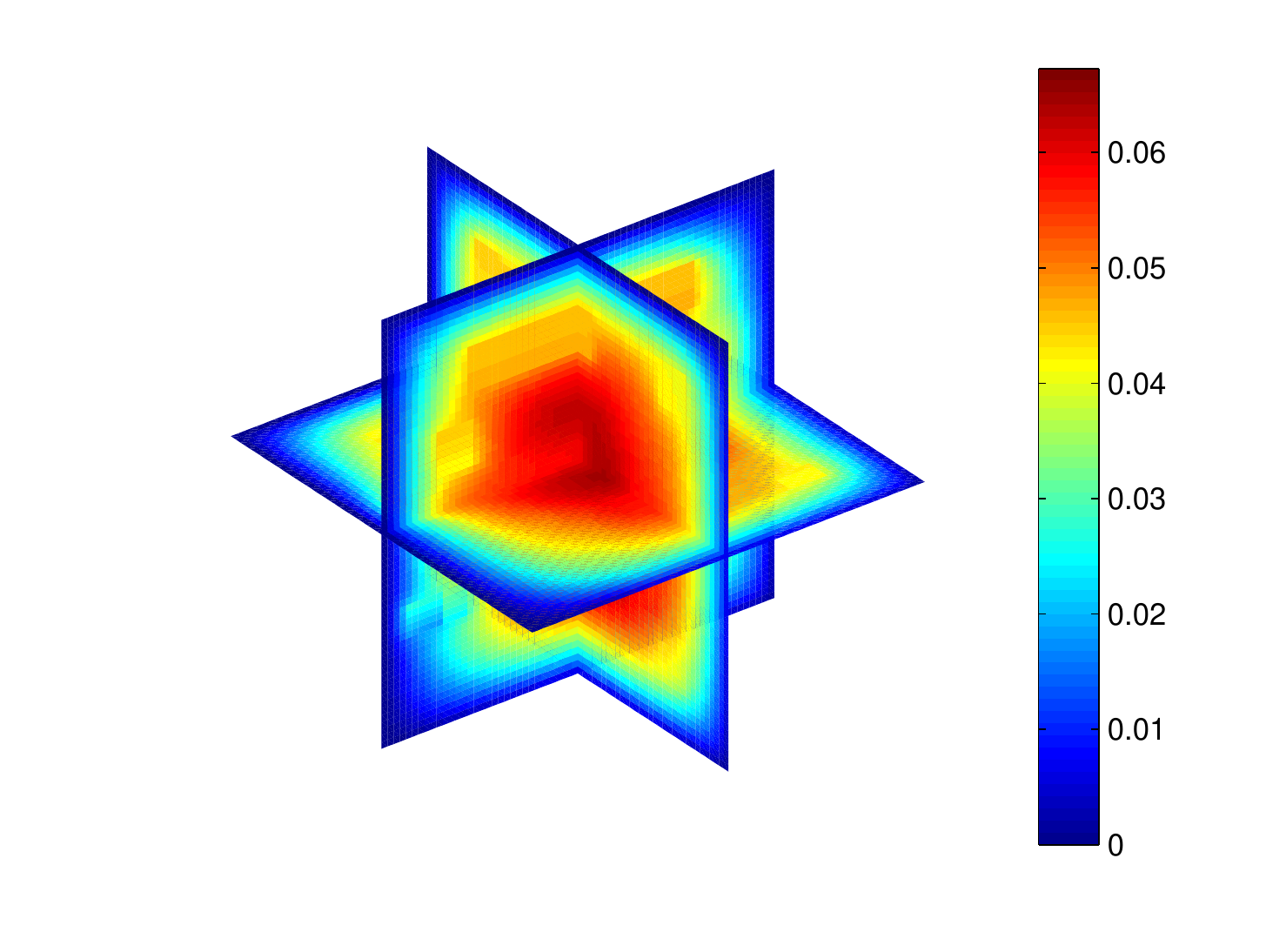}}
	\subfigure[Third component]{
		\includegraphics[width=2in]{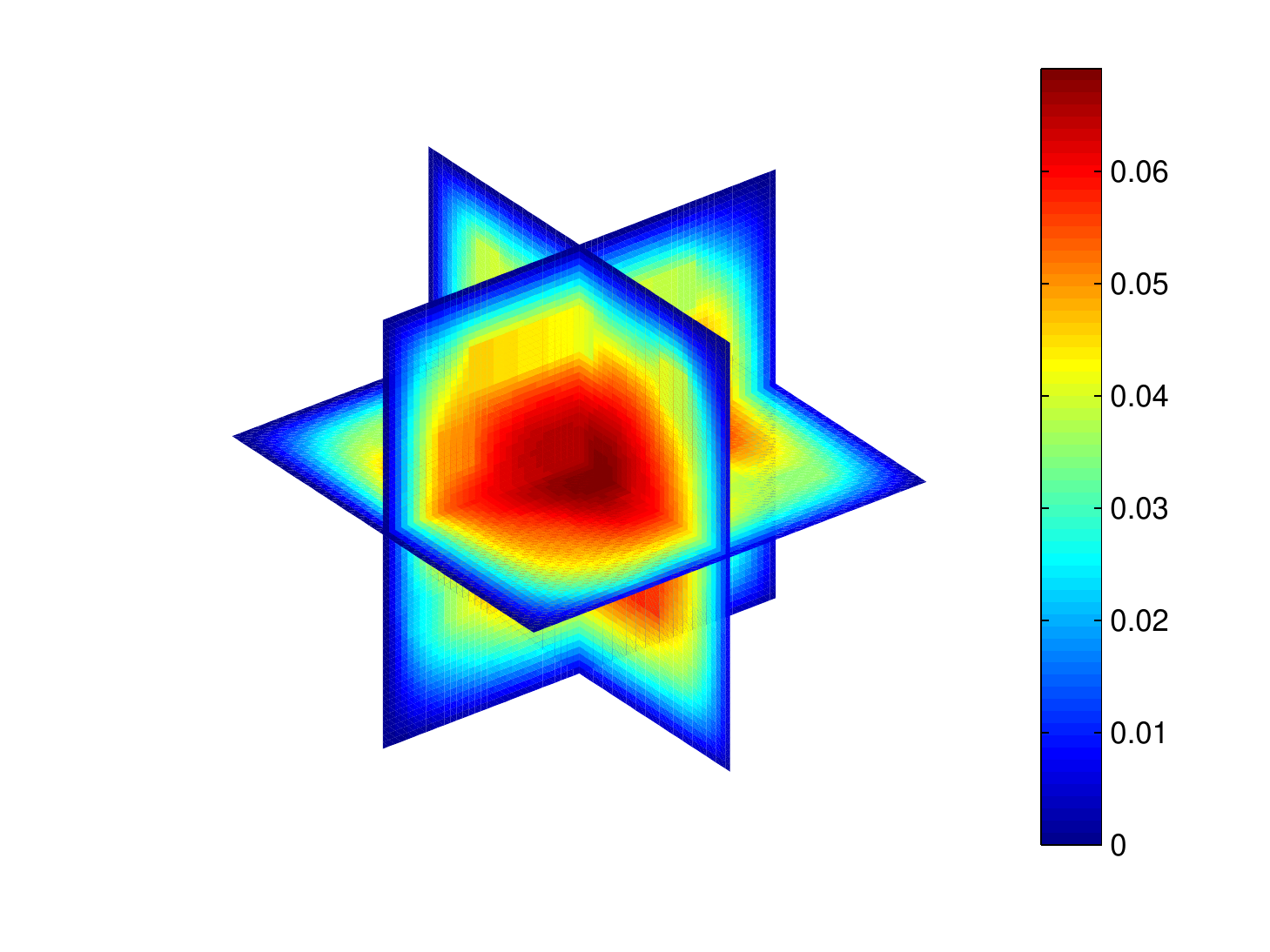}}	
	\caption{Reference solution}
	\label{fig:finedis3d}
\end{figure}

\begin{figure}[H]
	\centering
	\subfigure[First component]{
		\includegraphics[width=2in]{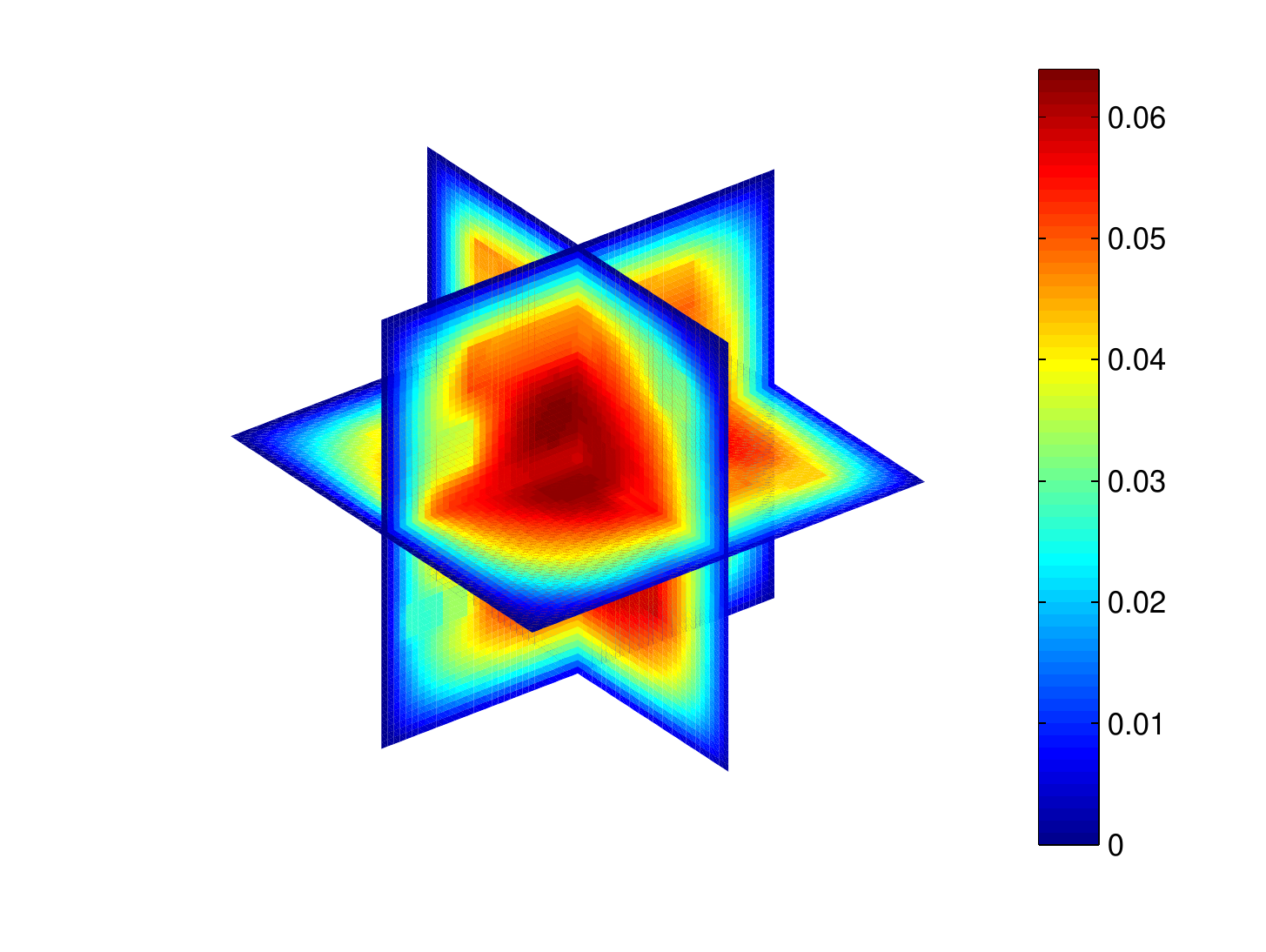}}
	\subfigure[First component]{
		\includegraphics[width=2in]{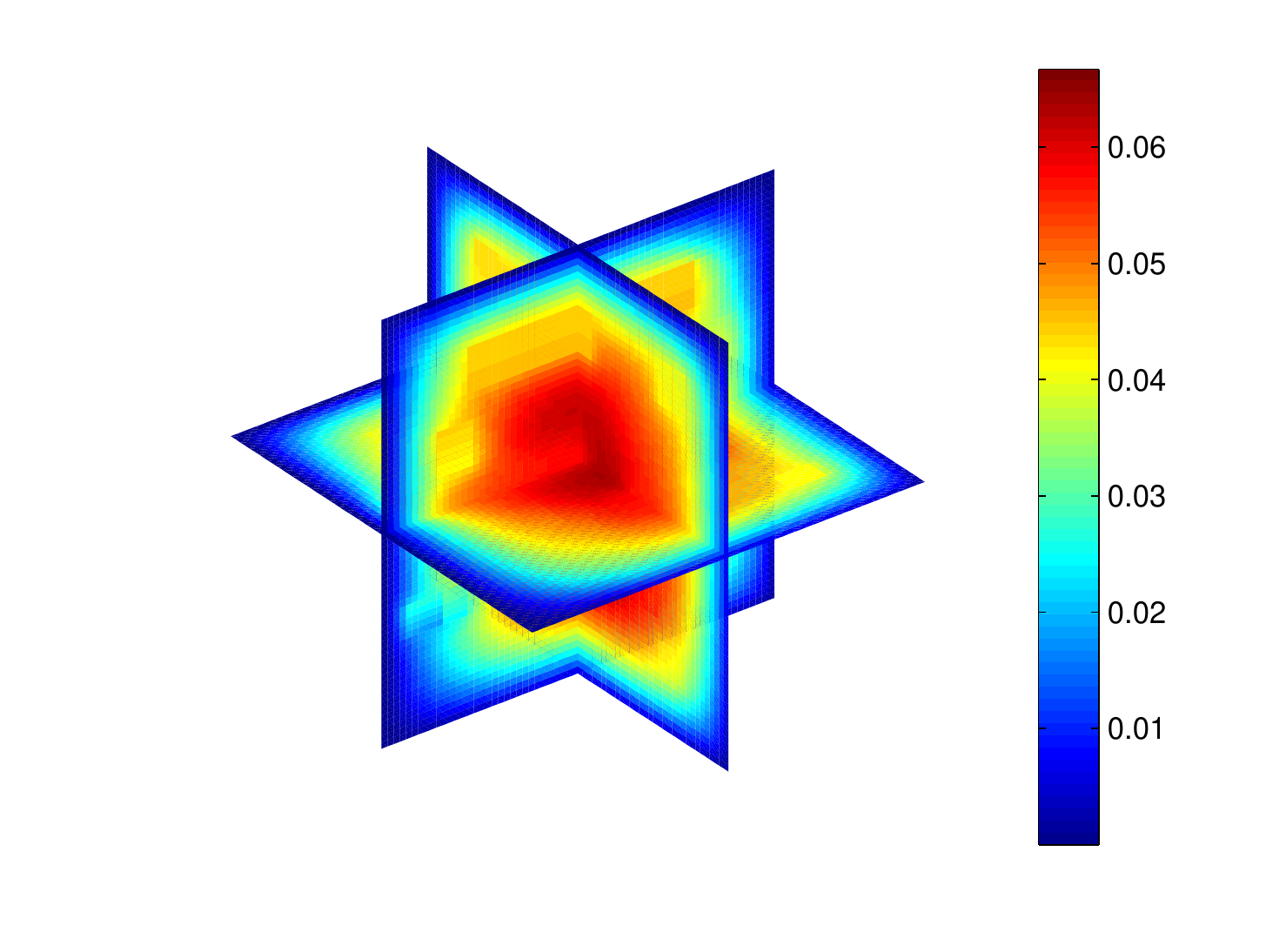}}
	\subfigure[Third component]{
		\includegraphics[width=2in]{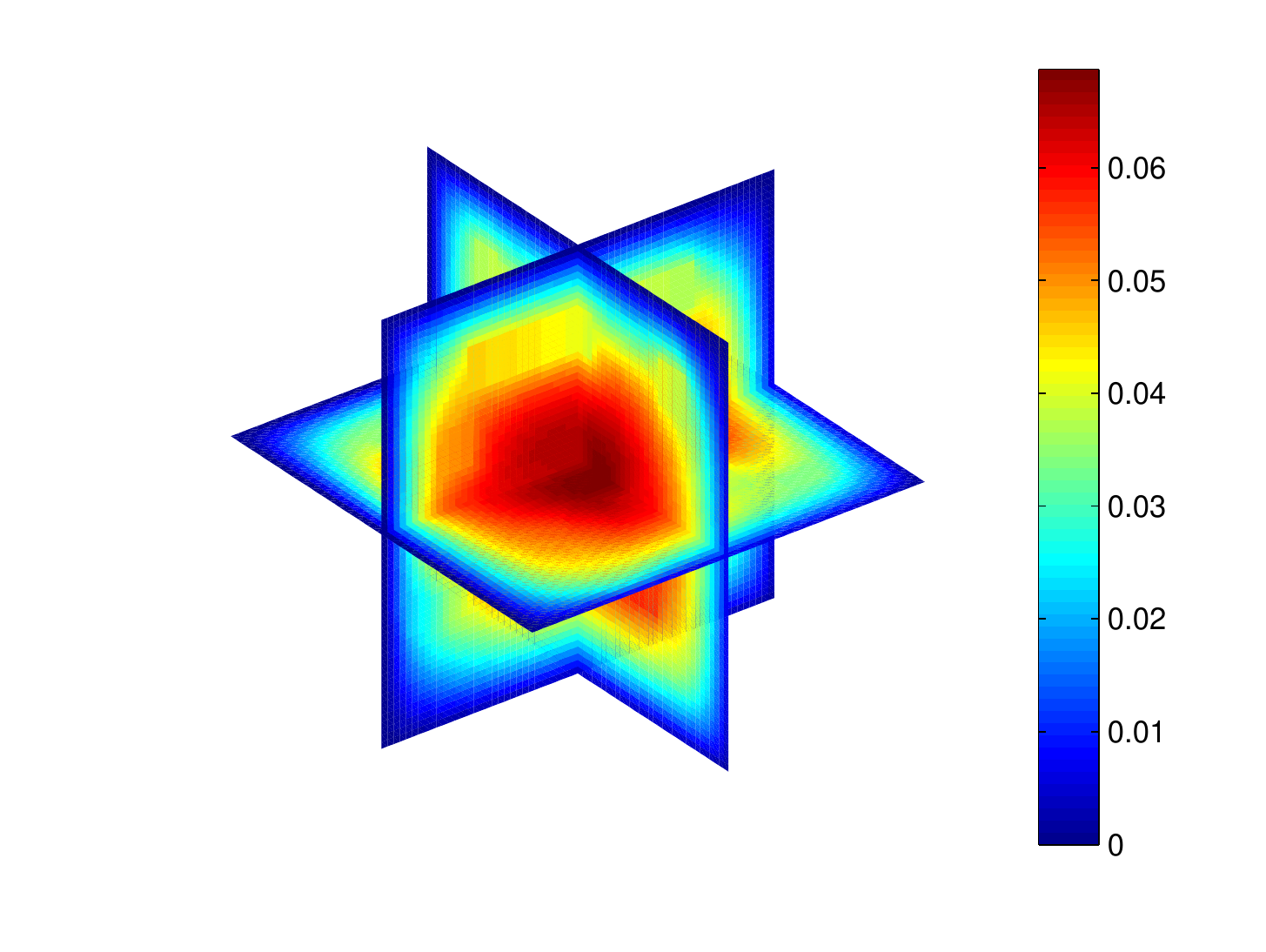}}	
	\caption{Relaxed CEM-GMsFEM solution, $N_b$=8, H$=1/16$, $n_\text{ov}=3$.}
	\label{fig:msdis3d} 
\end{figure}

\begin{figure}[H]
	\centering
	\includegraphics[width=2.5in]{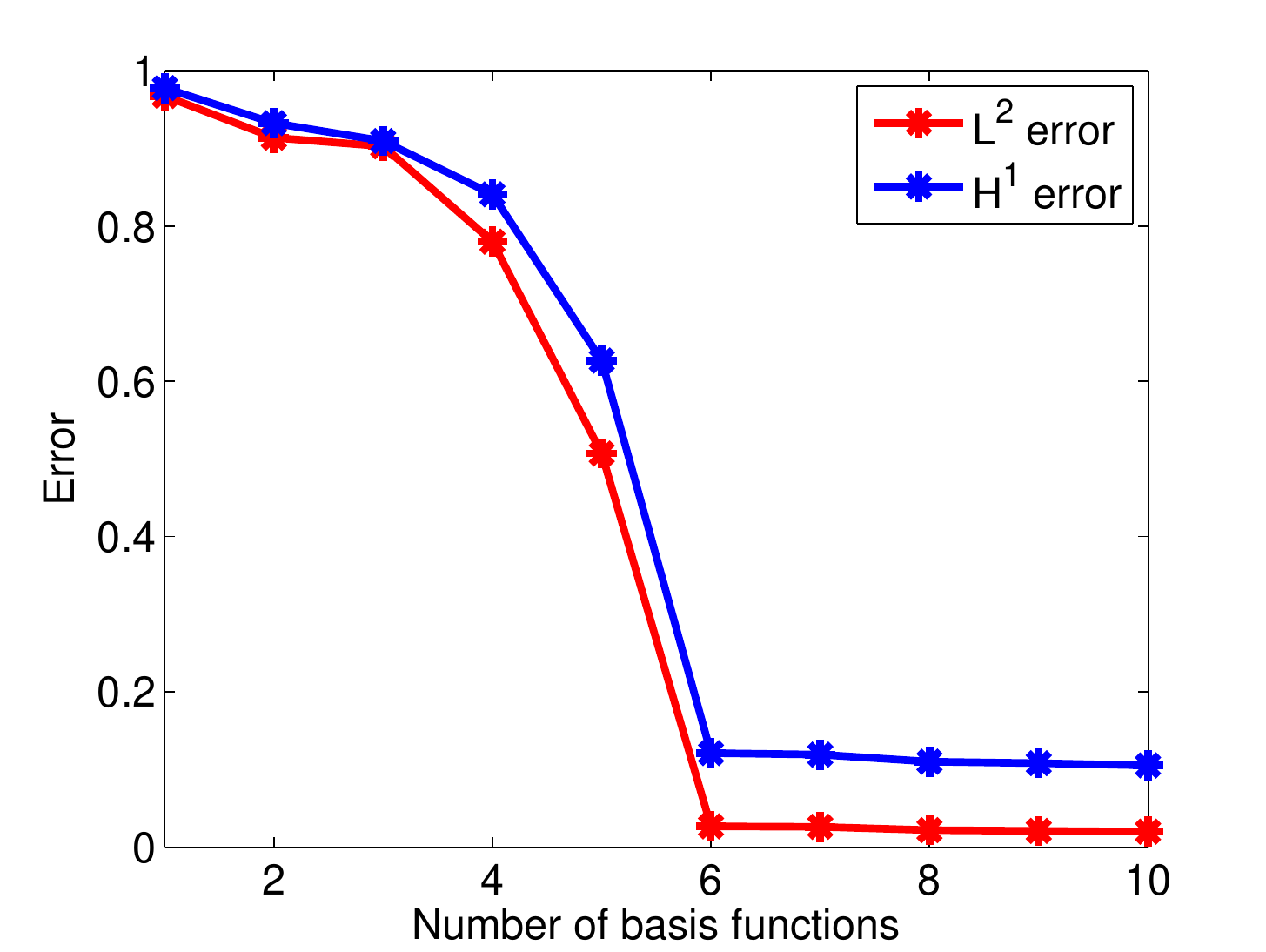}
	\caption{Numerical results  with different numbers of basis functions,  $H=1/16$, $n_\text{ov}=3$, relaxed case.}
	\label{fig:vb3d}
\end{figure}

In subsection, we present the test results of relaxed CEM-GMsFEM on model 2. 
We consider using a oversampling layer of $2\text{log}(H)/\text{log}(1/8)$ and 8 eigenfunction in local
auxiliary basis space. The results with varying coarse grid size and fixed $N_b$  are shown in Table \ref{ta:3dH2}. We also observe that multiscale solution converges with respect to the coarse grid size.
Figure \ref{fig:finedis3d} and Figure \ref{fig:msdis3d} show the displacement fields comparison between
the reference solution and multiscale solution, we can see that multiscale solution can approximate the 
reference solution pretty well. 
We also consider using different number of auxiliary basis functions, the results are shown in Figure \ref{fig:vb3d}. Again we find that using more basis
will increase the accuracy of the coarse solution. Once the basis number
reaches a value, the decay of the error becomes slower.
By varying the contrast of the media, we obtain various results shown in Table \ref{ta:3dcontrast}. We  find that increasing the number of oversampling layers can increase robustness of the method.
The uniform and adaptive online convergence history are presented in Table \ref{ta:3donline} and \ref{ta:3donline1} respectively, we can observe a fast decay of the error with more basis used, adaptive enrichment is better than
the uniform enrichment in reducing the $L^2$ error. The reason that why uniform is better in reducing 
$H^1$ error is that using uniform number of basis functions may yield the more  smoother solution than non-uniform case, the relative errors are more uniform in different regions. 
From the results shown in Table \ref{ta:3donline1} and Figure \ref{fig:vb3d}, we can observe the 
superiority of the online basis.

\begin{table}[H]
\centering \begin{tabular}{|c|c|c|c|c|c|}\hline
$N_b$ &$H$& $n_\text{ov}$  & $e_{L^2}$   & $e_{H^1}$   \tabularnewline\hline
8&1/8	&2&1.73e-01&3.08e-01 \tabularnewline\hline
8&1/16	&3&2.15e-02&1.10e-01 \tabularnewline\hline
\end{tabular}
\caption{Numerical results with varying coarse grid size $H$ for the test model 2, relaxed case.}
\label{ta:3dH2}
\end{table}

\begin{table}[H]
\centering
\begin{adjustbox}{max width=\textwidth}		
\begin{tabular}{|c|c|c|c|c|c|}\hline
\diagbox{$n_\text{ov}$}{$E_2$} &$10^2$ &$10^4$ &$10^6$  \tabularnewline\hline
3&2.59e-02&1.10e-01&5.64e-01 \tabularnewline\hline
4&2.21e-02&2.66e-02&1.28e-01 \tabularnewline\hline
\end{tabular}
\end{adjustbox}
	\caption{Comparison of ($e_{H^1}$) various number of oversampling layers and different contrast values for test model 2, relaxed case, $H=1/16, N_b=8$.}
	\label{ta:3dcontrast} 
\end{table}

\begin{table}[H]
	\centering \begin{tabular}{|c|c|c|c|c|c|}\hline
		Dof &H&$n_\text{ov}$ & $e_{L^2}$   & $e_{H^1}$   \tabularnewline\hline
		24576&1/16	&3&2.63e-02&1.21e-01     \tabularnewline\hline
		27951&1/16	&3&8.61e-05 &4.04e-03 \tabularnewline\hline
		31326&1/16	&3&1.94e-06&1.02e-03   \tabularnewline\hline
	\end{tabular}
	\caption{Uniform enrichment error decay history for the test model 2, $H=1/16$, 6 offline basis used.}
	\label{ta:3donline} 
\end{table}

\begin{table}[H]
	\centering \begin{tabular}{|c|c|c|c|c|c|}\hline
		Dof &H&$n_\text{ov}$ & $e_{L^2}$   & $e_{H^1}$   \tabularnewline\hline
		24576&1/16	&3&2.63e-02&1.21e-01     \tabularnewline\hline
		25521&1/16	&3&2.43e-04&2.54e-02\tabularnewline\hline
		26602&1/16	&3&4.51e-05&1.31e-02   \tabularnewline\hline
		27798&1/16	&3&3.10e-05&1.02e-02   \tabularnewline\hline
	\end{tabular}
	\caption{Adaptive enrichment with $\theta=0.1$ error decay history for the test model 2, $H=1/16$, 6 offline basis used.}
	\label{ta:3donline1} 
\end{table}

\section{Conclusions}\label{sec:conclusion}
In this paper, we propose Constraint Energy Minimizing GMsFEM for solving linear elasticity problems in high-contrast media. We introduce the construction of offline and online Constraint Energy Minimizing
multiscale basis functions. To construct the offline basis, we first construct an auxiliary space, and then
solve energy minimizing problems in target coarse block. The online basis is construct via solving a local 
problem in a oversampling domain with the residual as source.  We provided numerical tests on 2D and 3D models 
to demonstrate the accuracy of our method.

\section*{Acknowledgements}

EC's work is partially supported by Hong Kong RGC General Research Fund (Project 14304217) and CUHK Direct Grant for Research 2017-18.

\bibliographystyle{plain}
\bibliography{references}

\end{document}